\documentclass[12pt,reqno]{amsart}
\usepackage{amsmath,amsfonts,amsthm,amsopn,amssymb}
\usepackage{cite,marginnote}
\usepackage{color,enumitem,graphicx}
\usepackage[colorlinks=true,urlcolor=blue,
citecolor=red,linkcolor=blue,linktocpage,pdfpagelabels,
bookmarksnumbered,bookmarksopen]{hyperref}
\usepackage[english]{babel}

\usepackage[left=2.9cm,right=2.9cm,top=2.8cm,bottom=2.8cm]{geometry}
\usepackage[hyperpageref]{backref}

\numberwithin{equation}{section}

\makeindex
\newtheorem{theorem}{Theorem}[section]
\newtheorem{lemma}{Lemma}[section]
\newtheorem{corollary}{Corollary}[section]
\newtheorem{proposition}{Proposition}[section]
\newtheorem{remark}{Remark}[section]
\newtheorem{example}{Example}[section]
\newtheorem{definition}{Definition}[section]

\newtheorem{problem}{Problem}

\newcommand{\iy}{\infty}

\newcommand{\s}{\section}

\newcommand{\DD}{\Delta}

\newcommand{\la}{\lambda}

\newcommand{\R}{\mathbb R}
\newcommand{\al}{\alpha}

\newcommand{\rg}{\rightarrow}

\newcommand{\bt}{\begin{theorem}}
\newcommand{\et}{\end{theorem}}
\newcommand{\bl}{\begin{lemma}}
\newcommand{\el}{\end{lemma}}
\newcommand{\bd}{\begin{definition}}
\newcommand{\ed}{\end{definition}}
\newcommand{\bc}{\begin{corollary}}
\newcommand{\ec}{\end{corollary}}
\newcommand{\bp}{\begin{proof}}
\newcommand{\ep}{\end{proof}}
\newcommand{\bx}{\begin{example}}
\newcommand{\ex}{\end{example}}
\newcommand{\bi}{\begin{exercise}}
\newcommand{\ei}{\end{exercise}}
\newcommand{\bo}{\begin{proposition}}
\newcommand{\eo}{\end{proposition}}
\newcommand{\br}{\begin{remark}}
\newcommand{\er}{\end{remark}}
\newcommand{\be}{\begin{equation}}
\newcommand{\ee}{\end{equation}}
\newcommand{\ba}{\begin{align}}
\newcommand{\ea}{\end{align}}
\newcommand{\bn}{\begin{enumerate}}
\newcommand{\en}{\end{enumerate}}
\newcommand{\bg}{\begin{align*}}
\newcommand{\bcs}{\begin{cases}}
\newcommand{\ecs}{\end{cases}}

\newcommand{\bean}{\begin{eqnarray*}}
\newcommand{\eean}{\end{eqnarray*}}

\newcommand{\N}{\mathcal {N}}
\def\Proof{\noindent{\bf Proof}\quad}
\def\qed{\hfill$\square$\smallskip}

\title[A perturbation approach to studying sign-changing solutions]{A perturbation approach to studying sign-changing solutions of Kirchhoff equations with a general nonlinearity}

\author[Z.~Liu]{Zhisu Liu}
\author[Y.~Lou]{Yijun Lou}
\author[J.~Zhang]{Jianjun Zhang}

\address[Z.\ Liu]{\newline\indent School of Mathematics and Physics
\newline\indent University of South China
\newline\indent 421001, Hengyang, Hunan, P. R. China
}
\email{\href{mailto:liuzhisu183@sina.com}{liuzhisu183@sina.com}}

\address[Y.\ Lou]{\newline\indent Department of Applied Mathematics
\newline\indent The Hong Kong Polytechnic University
\newline\indent Hung Hom, Kowloon, Hong Kong
}
\email{\href{mailto:yijun.lou@polyu.edu.hk}{yijun.lou@polyu.edu.hk}}

\address[J.~Zhang]{\newline\indent College of Mathematics and Statistics
\newline\indent
Chongqing Jiaotong University
\newline\indent
Chongqing 400074, PR China
\newline\indent and
\newline\indent Dip. di Scienza e Alta Tecnologia
\newline\indent
Universit\`{a} degli Studi dell'Insubria
\newline\indent
via G.B. Vico 46, 21100 Varese, Italy}
\email{\href{mailto:zhangjianjun09@tsinghua.org.cn}{zhangjianjun09@tsinghua.org.cn}}

\thanks{Z. Liu was partially supported by the NSFC (Grant No.11271115),
NSF of Hunan Province (No. 2017JJ3265). J. Zhang is the corresponding author and was supported by the NSFC (Grant No. 11871123).}


\begin{document}

\begin{abstract}
By employing a novel perturbation approach and the method of invariant sets of descending flow, this manuscript investigates the existence and multiplicity of sign-changing solutions to a class of semilinear Kirchhoff equations in the following form
$$
-\left(a+ b\int_{\R^3}|\nabla u|^2\right)\triangle {u}+V(x)u=f(u),\,\,x\in\R^3,
$$
where $a,b>0$ are constants, $V\in C(\R^3,\R)$, $f\in C(\R,\R)$. The methodology proposed in the current paper is  robust, in the sense that,
the monotonicity condition for the nonlinearity $f$ and the coercivity condition of $V$ are not required. Our result improves the study made by Y. Deng, S. Peng and W. Shuai ({\it J. Functional Analysis}, 3500-3527(2015)), in the sense that, in the present paper, the nonlinearities include the power-type case $f(u)=|u|^{p-2}u$ for $p\in(2,4)$, in which case, it remains open in the existing literature that whether there exist infinitely many sign-changing solutions to the problem above without the coercivity condition of $V$. Moreover, {\it energy doubling} is established, i.e., the energy of sign-changing solutions
is strictly large than two times that of the ground state solutions for small $b>0$.
\end{abstract}

\maketitle

\s{Introduction}
\renewcommand{\theequation}{1.\arabic{equation}}
\subsection{Background}
 In the present paper, we investigate the sign-changing solutions of the following Kirchhoff equation
\begin{align}\label{K}\tag{K}
-\left(a+ b\int_{\R^3}|\nabla u|^2\right)\triangle {u}+V(x)u=f(u),\,\,x\in\R^3,\,\,u\in H^1(\R^3),
\end{align}
where $V\in C(\R^3,\R)$, $f\in C(\R,\R)$, and $a,b>0$ are positive constants.
Problem (K) arises in an interesting physical context.
Indeed, if we set $V(x)=0$ and a domain $\Omega\subset\R^3$ and replace $f(u)$ by $f(x,u)$,
problem (K) becomes as
the following Dirichlet problem:
\begin{equation} \label{eqn:youjielv}
\left\{
  \begin{aligned}
    -\left(a+b\int_{\Omega}|\nabla u|^2\right)\triangle {u} =f(x,u),\quad & \mbox{in}\,\,\Omega, \\
    u=0,\quad & \mbox{on}\,\, \partial\Omega,
   \end{aligned}
\right.
\end{equation}
which is the general form of the stationary counterpart of the hyperbolic Kirchhoff equation
\begin{equation}\label{eqn:stationary}
\rho\frac{\partial^2u}{\partial t^2}=\left[\frac{P_0}{h}+\frac{E}{2L}\int_0^L\left(\frac{\partial u}{\partial x}\right)^2dx\right]\frac{\partial^2u}{\partial x^2}.
\end{equation}
Such equations were proposed by Kirchhoff in \cite{Kirchhoff83} as an
existence of the classical D'Alembert's wave equations for free
vibration of elastic strings with taking into account
the changes in length of the string produced by transverse
vibrations. In (\ref{eqn:stationary}), $u$
denotes the displacement, $b$ is the initial tension while $a$ is
related to the intrinsic properties of the string (such as Young¡¯s
modulus). The nonlinearity $f(x,u)$ stands for the external force. Besides,
we also point out that Kirchhoff problems appear in other fields
like biological systems, such as population density, where $u$
describes a process which depends on the average of itself.
For the further physical background, we refer the readers to \cite{Cavalcanti01,Chipot97}.

\subsection{Overview and motivation}
Due to the presence of the term $\int_{\Omega}|\nabla u|^2$, equations (\ref{eqn:youjielv}) and (\ref{eqn:stationary}) are no longer
a pointwise identity and therefore, Kirchhoff problems are viewed as being nonlocal.
This observation brings mathematical challenges to the analysis, and at the same time,
makes the study of such a problem particularly interesting.  Another difficulty lies
on showing the boundedness of Palais-Smale sequences.  In the past decades, this kind
of problems have been receiving extensive attention. Initiated by Lions \cite{Lions78},
by using the variational methods, the solvability of Kirchhoff
type equation (\ref{eqn:youjielv}) has been investigated in many studies
(see \cite{Alves10,Alves05,Arosio96,Liang13,Ma03,Mao09,Perera06,Sun16,Zhang06} and the references therein).
There also have been many interesting works about the existence of positive solutions,
multiple solutions, ground states and semiclassical states to Kirchhoff type equation (K)
via variational methods, see for instance
\cite{Alves12,Azzollini11,He11,Li14,Figueiredo14,Liu15,Nie12,Wang12,Wu11,He15,He16,Xie17,Liu17} and the references therein.

In recent years, another interesting topic is the existence of sign-changing solutions
to Kirchhoff problems. Via the minimax approach and the method of invariant sets of descent flow,
Zhang and Perera \cite{Zhang06} and Mao and Zhang \cite{Mao09} proved the existence of sign-changing
 solutions of (\ref{eqn:youjielv}) provided that the function $f$ satisfies the $4$-superlinear
 growth condition:
\begin{align}\label{11111}\tag{4-superlinear}
\lim_{|t|\rightarrow+\infty}F(x,t)/t^4=+\infty,\,\mbox{uniformly for}\,\,x\in\Omega,
\end{align}
where $F(x,t)=\int_{0}^tf(x,s)ds$.
Meanwhile, the authors in \cite{Zhang06}  also considered the asymptotically 4-linear case:
\begin{align}\label{11112}\tag{asymptotically 4-linear}
\lim_{|t|\rightarrow+\infty}f(x,t)/bt^3=\kappa>0,\,\mbox{uniformly for}\,\,x\in\Omega
\end{align}
or the 4-sublinear case:
\begin{align}\label{11113}\tag{4-sublinear}
|f(x,t)|\le C(1+|t|^{p-1}),\, x\in\Omega,\, t\in\R,\,\mbox{for some}\,\,p\in(2,4).
\end{align}
Here we should point out that the 4-sublinear case in \cite{Zhang06} poses an additional restriction on $f$: $F(x,t)\ge c t^2$ for $|t|$ small, which rules out  cases including that $f(x,t)=|t|^{p-2}t$ for $p\in(2,4)$.

Subsequently, by the constraint variational method, Shuai \cite{Shuai15} obtanied the existence of least energy sign-changing solutions to problem (\ref{eqn:youjielv}). The author also showed that the energy of any sign-changing solutions is strictly larger than that of the ground state solutions of (\ref{eqn:youjielv}).  Note that one key assumption in \cite{Shuai15} is on the Nehari type monotonicity condition of $f$:
\begin{center}
{\bf (Ne)} $\frac{f(t)}{|t|^3}$ is increasing on $t\in(-\infty,0)\cup(0,+\infty)$.
\end{center}
Recently, a weaker condition of (Ne) was proposed by Tang and Cheng \cite{Tang16}, based on which the authors established the existence of least energy
sign-changing solution to problem (\ref{eqn:youjielv}) and showed that the energy of any
sign-changing solutions is strictly  two times larger than  that of the ground state solutions.
Recently, by virtue of the invariant sets method, Sun et al. \cite{Sun18} obtained infinitely many sign-changing solutions of problem (K) when
$f(t)\thicksim|t|^{p-2}t$ as $t\rg\iy$ for $p\in(2,4]$. However, $V$ is required to be coercive,
$\lim_{|x|\rightarrow\infty}V(x)=\infty$, which plays a key role in obtaining the boundedness of Palais-Smale equences. Without the coercive condition, Deng, Peng and Shuai \cite{Deng14} established the existence and asymptotic behavior of nodal solutions to (K) with the nonlinear term $f(|x|,t)$. Precisely, they obtained the existence of a sign-changing solution, which changes signs exactly $k$ times for any $k\in \mathbb{N}$. Their procedure of arguments is to transform the original problem to solving a system of $(k+1)$ equations with $(k+1)$ unknown functions
$u_i$ with disjoint supports. Then the nodal solution is constructed through gluing $u_i$ by matching the normal derivative at each junction point. We highlight that the monotonicity condition like (Ne) plays a crucial role in \cite{Deng14}. Moreover, they further assumed that $F(|x|,u)/u^4\rg+\iy$ as $|u|\rg\iy$. We should point our that the case $f(t)=|t|^{p-2}t$ for $p\in(2,4)$ remains open in \cite{Deng14}. Very recently, this case was considered in \cite{Cassani19} for a class of Kirchhoff type equation $\left[1+\la\int_{\R^3}|\nabla u|^2+V(x)u^2\right]\left[-\DD u+V(x)u\right]=f(u)$ in $R^3$ and the coercive condition was not imposed. When the coefficient $\la$ is sufficiently small, multiple sign-changing solutions were obtained.

\subsection{Open question} The results \cite{Sun18,Cassani19,Deng14} suggest the following open question.
\begin{problem}
Does problem (K) admit sign-changing solutions of (K) without (Ne) or the coercivity condition of $V$? In particular, without the coercivity condition, does there exist infinitely many sign-changing solutions of (K) in the case $f(t)=|t|^{p-2}t$ for $p\in(2,4)$?
\end{problem}

\section{Main results} \label{sec2}
\renewcommand{\theequation}{2.\arabic{equation}}
The main interest of the present paper is to give an affirmative answer to this open question.
\vskip0.1in
\subsection{Variational setting} Throughout this paper, we assume the external Schr\"odinger potential $V\in C(\R^3,\R)$ enjoys the following condition
\begin{itemize}
\item[ ($V_1$)] $V(x)=V(|x|)$ for any $x\in\R^3$, and $\inf\limits_{x\in\R^3}V(x):=V_0>0$,
\end{itemize}
and $f\in C(\R,\R)$ satisfies the following hypotheses
\begin{itemize}
\item[ ($f_1$)]$f\in C(\R,\R)$ and $\lim\limits_{t\rightarrow0}\frac{f(t)}{t}=0$;
\item[ ($f_2$)]$\limsup\limits_{|t|\rightarrow\infty}\frac{|f(t)|}{|t|^{p-1}}<\infty$ for some $p\in(2,6)$;
\item[ ($f_3$)]there exists $\mu>2$ such that $tf(t)\geq\mu F(t)>0$ for $t\not=0$, where $F(t)=\int_{0}^tf(s)ds$.
\end{itemize}
\begin{remark}
By ($f_2$) and ($f_3$), we know $2<\mu\leq p<6$. As a reference model, $f(u)=|u|^{p-2}u$ satisfies ($f_1$)-($f_3$) for $p\in(2,6)$.
\end{remark}

To proceed, we first define the Hilbert space
\begin{align*}
 E=
\left\{u\in H_r^1(\R^3):\int_{\R^3} V(x)u^2<\infty\right\}
\end{align*}
with the inner product
$$
\langle u,v\rangle=\int_{\R^3}a\nabla
u\nabla v+V(x)uv
$$
and the norm
$$
\|u\|:=\sqrt{\langle u,u\rangle}=\left(\int_{\R^3}a|\nabla
u|^2+V(x)u^2\right)^{\frac{1}{2}}.
$$
Obviously, it follows from ($V_1$) that the embedding $E\hookrightarrow L^q(\R^3)$ is compact for $2<q<6$ (see Strauss \cite{Strauss77}). The associated energy functional $I:E\rightarrow\R$ is given by
\begin{equation*}\label{eqn:fanhan}
I(u)=\frac{1}{2}\|u\|^2
+\frac{b}{4}\left(\int_{\R^3}|\nabla u|^2\right)^2-\int_{\R^3} F(u).
\end{equation*}
It is a well-defined $C^1$ functional in $E$ and its derivative is given by
\begin{equation*}\label{eqn:Frechet}
 I'(u)v=\langle u,v\rangle+b\int_{\R^3}|\nabla u|^2\int_{\R^3}\nabla u\nabla v -\int_{\R^3} f(u)v,\,\,v\in E.
\end{equation*}
\begin{definition}
$u$ is called a weak solution of (K), if $u\in E$ satisfies $ I'(u)\varphi=0$ for all $\varphi\in C_0^\infty(\R^3)$. Furthermore, if $u^\pm\not\equiv0$, then $u$ is called a sign-changing solution of (K), where
$u^\pm=\max\{\pm u,0\}$.
\end{definition}

Before stating our main results, we impose some additional hypotheses on $V(x)$ as follows.
\begin{itemize}
\item[($V_2$)] $V(x)\in C(\R^3,\R^+)$ is differentiable and satisfies $(\nabla V(x), x)\in L^{\infty}(\R^3)\cup L^{\frac{3}{2}}(\R^3)$. Moreover,  there exists $\mu>2$ such that
$$
\frac{\mu-2}{\mu}V(x)-(\nabla V(x), x)\geq0;
$$
\end{itemize}

Under these assumptions on $V(x)$ and $f(x)$, we establish main results in terms of existence, multiplicity and
energy doubling on solutions, which are stated in the following three subsections.

\subsection{Existence}

Our first result reads as follows.
\begin{theorem}{\bf(Existence)}\label{Thm:existence}
If ($V_1$)-($V_2$) and ($f_1$)-($f_3$) hold, then problem (K)
admits at least one radially symmetric ground state sign-changing solution.
\end{theorem}

Observe that if $b=0$, problem (K) reduces to the following {\it local} Schr\"odinger equation
\begin{equation}\label{eqn:ellptic}
-a\triangle {u}+V(x)u=f(u),
\end{equation}
which does not depend on the nonlocal term $\int_{\R^3}|\nabla u|^2$ any more. To look for sign-changing solutions of equation (\ref{eqn:ellptic}), we list several approaches introduced in the literature. Based on the Nehari manifold technique,
Cerami, Solimini and Struwe \cite{Cerami86} proved the existence of sign-changing
solutions for elliptic problems involving critical
exponent(see also \cite{Bartsch01,Cao88}). The heat flow method was explored in
\cite{Chang03} to study the existence of sign-changing solutions. Morse theory can also be used
to consider the existence of sign-changing solutions (see \cite{Chang04}). In finding sign-changing
solutions of elliptic problems, the method of invariant sets of descending flow has been a powerful tool. Here we refer to \cite{Bartsch04,Bartsch05,Zou} and the reference therein. However, we should address a remark on the case $b=0$.

\begin{remark} However, in contrast to problem (\ref{eqn:ellptic}), the non-locality leads problem (K) to be more complicated in seeking sign-changing solutions.
\end{remark}

Now, we summarize the main difficulties and novelties developed in this study.
\begin{itemize}
\item[(1)] In finding sign-changing solutions of (\ref{eqn:ellptic}), a crucial ingredient is the following decomposition: for any $u\in E$,
\begin{equation}\label{eqn:fenjie1}
I_0(u)=I_0(u^+)+I_0(u^-),\,\, \langle I'_0(u), u^\pm\rangle=\langle I'_0(u^\pm), u^\pm\rangle,
\end{equation}
where $I_0$ is the energy functional of (\ref{eqn:ellptic}) defined by
$$
I_0(u)=\frac{a}{2}\int_{\R^3}|\nabla u|^2+\frac{1}{2}\int_{\R^3}V(x) u^2
-\int_{\R^3} F(u).
$$
However, due to the nonlocal term $\int_{\R^3}|\nabla u|^2$, one can get
\begin{equation}\label{eqn:fenjie3}
\left\{
\begin{array}{ll}
I(u)=I(u^+)+I(u^-)
+\frac{b}{2}\int_{\R^3}|\nabla u^+|^2\int_{\R^3}|\nabla u^-|^2, \vspace{0.2cm} \\
\langle I'(u), u^\pm\rangle=\langle I'(u^\pm), u^\pm\rangle+b\int_{\R^3}|\nabla u^+|^2\int_{\R^3}|\nabla u^-|^2,
\end{array}
\right.
\end{equation}
which do no longer satisfy the decomposition (\ref{eqn:fenjie1}). Motivated by \cite{Liuj15,Liuz16}, we attempt to find sign-changing solutions for problem (K) by using the method of invariants sets of a descending flow.

\item[(2)] For the case $f(t)=|t|^{p-2}t$,  $p\in(2,4]$, the effect of nonlocal term $\int_{\R^3}|\nabla u|^2$ results in two difficulties.
First, it seems much more complicated to find a similar auxiliary operator A (see \cite[Section 4]{Liuz16}),
 which plays a crucial role in constructing invariants sets of a descending flow associated
 with problem (K). A similar difficulty also arises in seeking sign-changing solutions of
 the Schr\"odinger-Possion systems
\begin{equation}\label{SP} \left\{
\begin{array}{ll}
-\Delta u+V(x)u+\phi u=|u|^{p-2}u&\mbox{in}\ \R^3,\\
-\Delta\phi=u^2&\mbox{in}\ \R^3.
\end{array}
\right.
\end{equation}
where $p\in(3,4)$. In \cite{Liuz16}, the authors overcome this difficulty
for $p\in(3,4)$ by adding a higher order local nonlinear term and the coercive
condition $(V_3)$.
Recently, this approach in \cite{Liuz16} was also used to deal with the Kirchhoff equation, see \cite{Sun18}.
Second, the so-called 4-(AR) condition fails, which makes
tough to get the boundedness of (PS) sequences. In \cite{Liuz16,Sun18}, the authors
recovered such boundedness due to the coercivity of $V(x)$.
However, without the coercivity condition,
the method in \cite{Liuz16,Sun18} is inapplicable any more. In this paper,
we develop particularly \emph{a new perturbation approach} by adding another perturbation,
which is nonlocal. For the perturbed problems, by minimax arguments in the presence of invariant
sets, we obtain sign-changing solutions. By passing to the limit, a convergence argument allows us to get sign-changing solutions of the original problem
(K), which involves the case $f(t)=|t|^{p-2}t$,
 $p\in(2,4)$. We emphasize that, without the coercive condition of $V$, this new
perturbation approach can also deal with system (\ref{SP}) in the case $p\in(3,4)$. Moreover,
we believe that this new perturbation approach should be of independent interest in other problems with the difficulty in verifying the boundedness of Palais-Smale sequences.

\end{itemize}

\subsection{Multiplicity}

Another aim of the paper is to prove the existence of infinitely many
sign-changing solutions to problem (K) when $f$ is odd.
\begin{theorem}{\bf(Multiplicity)}\label{Thm:many}
If (V$_1$)-($V_2$) and ($f_1$)-($f_3$) hold,
then problem (K) has infinitely many radially symmetric sign-changing solutions when $f$ is odd.
\end{theorem}

\br In proving Theorem \ref{Thm:many}, the main difficulties are three-fold. Firstly, since the nonlinearity is allowed to be $f(t)=|t|^{p-2}u$ with $p\in(2,4)$, it seems that the associated functional does not enjoy a linking structure. As a result, the minimax argument can not be used directly. Our strategy is that we adopt a perturbation approach by adding a higher order term $\beta|u|^{r-2}u$ to recover the linking structure. Secondly, since $\mu$ in $(f_3)$ may be smaller than $4$, without the coercive condition, the method in \cite{Liuz16} fails in proving the boundedness of Palais-Smale sequences. To overcome this obstacle, we give another perturbation term $\lambda\|u\|_2^{2\al}u$ in the left side of the equation. With the two perturbation terms and via the method of invariants sets of a descending flow, we obtain infinitely many sign-changing solutions $u_{\la,\beta}^k, k=1,2,\cdots$ of the perturbed problem as approximation solutions for the original equation. Then by passing to the limit, sign-changing solutions of the original problem are obtained. Lastly, since the minimax values $c_{\la,\beta}^k$ of the perturbed problem enjoy the different monotonicity properties on the two perturbation terms, it is not easy to distinguish the limits of $u_{\la,\beta}^k$ as $\la,\beta\rg0$. By using an auxiliary functional, we show that $c_{\la,\beta}^k\rg\iy$ as $k\rg\iy$ uniformly for $\la,\beta$. Based on this estimate, we obtain infinitely many radially symmetric sign-changing solutions for problem (K).
\er

\subsection{Energy doubling}

The last investigation is to establish {\it energy doubling} of sign-changing solutions
to problem (K) with $f(u)=|u|^{p-2}u$, $p\in(2,6)$. This fact has been proved for the {\it local} problem (\ref{eqn:ellptic}) when $V(x)$ is a constant or a periodic function \cite{Weth06}.  In particular, we denote the Nehari manifold associated with (\ref{eqn:ellptic}) by
$$
\mathcal{N}:=\{u\in E\setminus\{0\}:\,\,\langle I'_0(u),u\rangle=0\},
$$
and
\begin{equation} \label{eqn:p-ground}
c_0=\inf\{I_0(u):\,u\in\mathcal{N}\}.
\end{equation}
For any sign-changing solution $w\in E$ of (\ref{eqn:ellptic}),
it follows from the fact $w^\pm\in \mathcal{N}$ that
\begin{equation} \label{eqn:p-ground-}
I_0(w)=I_0(w^+)+I_0(w^-)\geq 2c_0.
\end{equation}
In fact, the minimizer of (\ref{eqn:p-ground}) is indeed a ground state solution of (\ref{eqn:ellptic}),
and $c_0>0$ is the ground state energy. If some sign-changing solution $w$ of (\ref{eqn:ellptic}) satisfies
$$
I_0(w)>2c_0,
$$
it was called in \cite{Weth11}(see also\cite{Weth06}) that $w$ satisfies \lq\lq energy doubling\rq\rq. They showed that any sign-changing solution of (\ref{eqn:ellptic}) satisfies energy doubling in the case $f(u)=|u|^{p-2}u$ for $p\in(2,6)$ and $V(x)$ is a constant or a periodic function \cite{Weth06}.

In the current study, we also estimate  the energy of sign-changing solutions to problem (K). Analogous to problem (\ref{eqn:ellptic}),  the definition of energy doubling  corresponding to problem (K) is given as follows.
\begin{definition}
Let $w_b\in E$ be a sign-changing solution of problem (K), we call $w_b$ satisfies energy doubling if $I(w_b)>2c_b$, where
\begin{equation} \label{eqn:p-ground1}
c_b=\inf\{I(u),u\in\mathcal{N}_b\} \text{ and } \mathcal{N}_b:=\{u\in E\setminus\{0\}:\,I'(u)=0\}.
\end{equation}
\end{definition}
Let $w_b\in E$ be a sign-changing solution of problem (K). Since the interaction of the positive and negative parts of solutions can not be neglected,
it follows from (\ref{eqn:fenjie3}) that
\begin{equation} \label{eqn:p-ground+}
w_b^{\pm}\not\in \mathcal{N}_b.
\end{equation}
Thus, a natural open question is whether energy doubling holds or not. Generally speaking, it is even not easy to compare $I(w_b)$ with $c_b$.
To proceed, we impose the additional assumptions on $V$.
\begin{itemize}
\item[(V$_4$)] $ V\in C^2([0,+\infty),\R^+)$ is radially symmetric and
$$
0<\inf_{r>0}V(r)\leq \sup_{r>0}V(r)<\infty;
$$
\item[(V$_5$)] $\inf_{r>0}[(V''(r)r^2+(3+\tau)V'(r)r+2\tau V(r)]>0$
\text{ with } $\tau=\frac{4(p-1)}{3+p}$.
\end{itemize}
By using an approximation procedure, we give a partial answer for such an open problem, that is, energy doubling holds if $b>0$ small.
Precisely, we have the following result.
\begin{theorem}{\bf(Energy doubling)}\label{Thm:ground}
For $f(u)=|u|^{p-2}u$, $p\in(2,6)$ and assume $V$ satisfies ($V_1$)-($V_2$) and ($V_4$)-(V$_5$),
then there exists $b^*>0$ such that, for any sign-changing solutions $w_b\in E$ of problem (K),
we have $I(w_b)>2c_b$ if $b<b^*$, i.e., energy doubling holds.
Furthermore, for any sequence $\{b_n\}$ with $b_n\rightarrow0$ as $n\rightarrow\infty$,
up to a subsequence, $w_{b_n}\rightarrow w_0$ in $E$, where $w_0$ is
a sign-changing solution of (\ref{eqn:ellptic}).
\end{theorem}
\begin{remark}
The assumptions (V$_4$)-(V$_5$) is imposed only to guarantee the uniqueness of positive solution
to equation (\ref{eqn:ellptic}) (see \cite{Tanaka99}), which is of use in our arguments.
\end{remark}

\subsection{Organization of this paper.}
The outline of our argument is as follows. A new perturbation approach is introduced, with which, we obtain the existence of sign-changing solutions to problem (K) in Section \ref{sec3}.  Section \ref{sec4} is devoted to proving Theorem \ref{Thm:many} by the minimax theorem through invariants sets of a descending flow. Finally, the energy doubling property is established for sign-changing solutions to problem (K)  in Section \ref{sec5}. The notations used in this paper are summarized as follows.
\vskip0.1in
\noindent{\bf Notations.}

\begin{itemize}
\item [$\bullet$] $\|u\|_p:=\big(\int_{\R^3}|u|^p\big)^{1/p}$ for $p\in [1,\infty)$.
\item [$\bullet$] $C$ will be used repeatedly to denote various
positive constants which may change from line to line.
\item [$\bullet$] $D^{1,2}(\R^{3}):=\left\{|\nabla u|\in L^2(\R^3):\,u\in L^6(\R^3)\right\}.$
\item [$\bullet$] $S$ denotes the best Sobolev constant, i. e.,
$$S:=\inf\limits_{u\in D^{1,2}(\R^{3})\setminus\{0\}}\frac{\int_{\R^3}|\nabla u|^2}{(\int_{\R^3}u^6)^{1/3}}.$$
\end{itemize}

\section{Existence}\setcounter{equation}{0}\label{sec3}
\renewcommand{\theequation}{3.\arabic{equation}}
In this section, we prove the existence of ground state sign-changing solutions to problem (K) when
$V(x)$ is a radial symmetric function.
\subsection{The perturbed problem}
Since we do not impose the well-known {\it Ambrosetti-Rabinowtiz} condition,
the boundedness of the Palais-Smale sequence becomes not easy to establish.
  A perturbed problem is introduced to overcome this difficulty.
Set $\alpha\in(0,\frac{\mu-2}{3\mu+2})$ and fix $\lambda$, $\beta\in(0,1]$ and
$r\in(\max\{p,\frac{9}{2}\},6)$, we consider the modified problem
$$
\left\{
  \begin{aligned}
-\left(a+ b\int_{\R^3}|\nabla u|^2\right)\triangle {u}+
V(x)u=f_{\lambda,\alpha,\beta}(u),\quad\mbox{in}\,\,\R^3,\\
u\in E,
  \end{aligned}
\right.\eqno{\rm{(K_{\lambda,\beta})}}
$$
where $$f_{\lambda,\alpha,\beta}(u)=f(u)+\beta|u|^{r-2}u-\lambda\left(\int_{\R^3}u^2\right)^{\alpha}u.$$ An associated functional can be constructed as
$$
I_{\lambda,\beta}(u)=I(u)+\frac{\lambda}{2(1+\alpha)}\left(\int_{\R^3}u^2\right)^{1+\alpha}-\frac{\beta}{r}\int_{\R^3}|u|^r.
$$
It is easy to show that $I_{\lambda,\beta}\in C^1(E,\R)$ and
$$
 I'_{\lambda,\beta}(u)v= I'(u)v +\lambda\left(\int_{\R^3}u^2\right)^{\alpha}\int_{\R^3}uv-\beta\int_{\R^3}|u|^{r-2}uv,\,\,u,v\in E.
$$

We will make use of the following Pohozaev type identity, whose proof is standard and
can be found in \cite{Berestycki1}.

\begin{lemma}\label{Lem:pohozaev}
Let $u$ be a critical point of $I_{\lambda,\beta}$ in $E$ for $(\lambda,\beta)\in
(0,1]\times(0,1]$, then
$$
\aligned
\frac{a}{2}\int_{\R^3}|\nabla u|^2+&\frac{3}{2}\int_{\R^N}V(x)u^2+\frac{1}{2}\int_{\R^3}(\nabla
V(x),x)u^2\\+&\frac{b}{2}\left(\int_{\R^3}|\nabla u|^2\right)^2
+\frac{3\lambda}{2}\left(\int_{\R^3}u^2\right)^{1+\alpha}-3\int_{\R^3}(F(u)+\beta|u|^r)=0.
\endaligned
$$
\end{lemma}
It is easy to prove that, for each $u\in E$, the following equation
\begin{equation} \label{eqn:in1-0}
    -\left(a+b\int_{\R^3}|\nabla u|^2\right)\triangle {v}+V(x)v+\lambda\left(\int_{\R^3}u^2\right)^{\alpha}v = f(u)+\beta|u|^{r-2}u
\end{equation}
has a unique weak solution $v\in E$.
In order to construct the descending flow for the functional $I_{\lambda,\beta}$, we introduce an auxiliary operator  $T_{\lambda,\beta}: u\in E\mapsto v\in E$,
where $v=T_{\lambda,\beta}(u)$ is the unique weak solution of problem (\ref{eqn:in1-0}).
Clearly, the fact that $u$ is a solution of problem (\ref{eqn:in1-0}) is equivalent to that $u$ is a fixed point
of $T_{\lambda,\beta}$, which is well defined based on the above arguments. Moreover, this operator is continuous, as stated in the next lemma.

\begin{lemma}\label{Lem:T1}
The operator $T_{\lambda,\beta}$ is well defined and continuous.
\end{lemma}
\Proof  Assume that $\{u_n\}\subset E$ with $u_n\rightarrow u$
strongly in $E$ as $n\rightarrow\infty$.
Let $v=T_{\lambda,\beta}(u)$ and
$v_n=T_{\lambda,\beta}(u_n)$, then we have
\begin{equation}\label{eqn:f-0}
\aligned
&\int_{\R^3}(a\nabla v_n\nabla w+V(x)v_nw)+b\int_{\R^3}|\nabla u_n|^2\int_{\R^3}\nabla v_n\nabla w
+\lambda\left(\int_{\R^3} u_n^2\right)^\alpha\int_{\R^3} v_nw
\\&=\int_{\R^3}f(u_n)w+\beta\int_{\R^3}|u_n|^{r-2}u_nw, \quad\forall w\in E
\endaligned
\end{equation}
and
\begin{equation}\label{eqn:f-00}
\aligned
&\int_{\R^3}(a\nabla v\nabla w+V(x)vw)+b\int_{\R^3}|\nabla u|^2\int_{\R^3}\nabla v\nabla w
+\lambda\left(\int_{\R^3} u^2\right)^\alpha\int_{\R^3} vw
\\&=\int_{\R^3}f(u)w+\beta\int_{\R^3}|u|^{r-2}uw, \quad\forall w\in E.
\endaligned
\end{equation}
We need to show $\|v_n-v\|\rightarrow0$ as $n\rightarrow\infty$. Indeed,
it follows from ($f_1$) and ($f_3$) that for any $\varepsilon>0$,
there exists $C_\varepsilon>0$ such that
\begin{equation}\label{eqn:f}
|f(t)|\leq \varepsilon |t|+C_\varepsilon|t|^{p-1}.
\end{equation}
Testing with $w=v_n$ in (\ref{eqn:f-0}) gives
$$
\aligned
&\|v_n\|^2+b\int_{\R^3}|\nabla u_n|^2\int_{\R^3}|\nabla v_n|^2+\lambda\|u_n\|_2^{2\alpha}
\int_{\R^3} v_n^2
\\&\leq\int_{\R^3}(\varepsilon |u_n|+C_\varepsilon|u_n|^{p-1})|v_n|
+\beta\int_{\R^3}|u_n|^{r-1}|v_n|,
\endaligned
$$
which, together with the H\"{o}lder inequality, imply that $\{v_n\}$ is a bounded sequence in $E$.
Assume $v_n\rightharpoonup v^*$ in $E$
and $v_n\rightarrow v^*$ in $L^s(\R^3)$ for $s\in(2,6)$ after extracting a subsequence, then by (\ref{eqn:f-0}) we have
\begin{equation}\label{eqn:f--}
\aligned
&\int_{\R^3}(a\nabla v^*\nabla w+V(x)v^*w)+b\int_{\R^3}|\nabla u|^2\int_{\R^3}\nabla v^*\nabla w
+\lambda\left(\int_{\R^3} u^2\right)^\alpha\int_{\R^3} v^*w
\\&=\int_{\R^3}f(u)w+\beta\int_{\R^3}|u|^{r-2}uw, \quad\forall w\in E.
\endaligned
\end{equation}
Hence $v^*$ is a weak solution of (\ref{eqn:in1-0}), which implies  $v=v^*$ by the uniqueness.
Moreover, taking $w=v_n-v$ in (\ref{eqn:f-0}) and (\ref{eqn:f-00})
and then subtracting, we have
\begin{equation}\label{eqn:f1}
\aligned
&\|v_n-v\|^2+b\int_{\R^3}|\nabla u_n|^2\int_{\R^3}|\nabla(v_n- v)|^2
+\lambda\|u_n\|^{2\alpha}\int_{\R^3} |v_n-v|^2\\
&=b\int_{\R^3}(|\nabla u_n|^2-|\nabla u|^2)\int_{\R^3}\nabla v\nabla(v_n-v)+
\lambda(\|u_n\|_2^{2\alpha}-\|u\|_2^{2\alpha})\int_{\R^3}v(v_n-v)\\
&+\int_{\R^3}(f(u_n)-f(u))(v_n-v)+\beta\int_{\R^3}(|u_n|^{r-2}u_n-|u|^{r-2}u)(v_n-v).
\endaligned
\end{equation}
It follows from (\ref{eqn:f--})-(\ref{eqn:f1}) and Sobolev's embedding inequality that
$v_n\rightarrow v$ in $E$ as $n\rightarrow\infty$. Therefore, $T_{\lambda,\beta}$ is continuous.
\qed

\begin{lemma}\label{Lem:T2}
\begin{itemize}
\item[\rm (1) ] $I'_{\lambda,\beta}(u)(u-T_{\lambda,\beta}(u))\geq \|u-T_{\lambda,\beta}(u)\|^2$ for all $u\in E$;
\item[\rm (2) ] $\|I'_{\lambda,\beta}(u)\|\leq \|u-T_{\lambda,\beta}(u)\|(1+C_1\|u\|^2+C_2\|u\|^{2\alpha})$
 for all $u\in E$, where $C_1$ and $C_2$ are two positive constants.
\end{itemize}
\end{lemma}
\Proof Since $T_{\lambda,\beta}(u)$ is the solution of equation (\ref{eqn:in1-0}), we have
$$
\aligned
I'_{\lambda,\beta}(u)(u-T_{\lambda,\beta}(u))=&\|u-T_{\lambda,\beta}(u)\|^2+b\int_{\R^3}|\nabla u|^2\int_{\R^3}|\nabla(u-T_{\lambda,\beta}(u))|^2\\
&+\lambda\|u\|_2^{2\alpha}\int_{\R^3}|u-T_{\lambda,\beta}(u)|^2,
\endaligned
$$
which means $I_{\lambda,\beta}'(u)(u-T_{\lambda,\beta}(u))\geq \|u-T_{\lambda,\beta}(u)\|^2$ for all $u\in E$.
Notice that for any $\varphi\in C_0^\infty(\R^3)$,
$$
\aligned
I_{\lambda,\beta}'(u)\varphi=&\int_{\R^3}[a\nabla(u-T_{\lambda,\beta}(u))\nabla\varphi
+V(x)(u-T_{\lambda,\beta}(u))\varphi]\\
&+b\int_{\R^3}|\nabla u|^2\int_{\R^3}\nabla(u-T_{\lambda,\beta}(u))\nabla \varphi
+\lambda\|u\|_2^{2\alpha}\int_{\R^3}(u-T_{\lambda,\beta}(u))\varphi,
\endaligned
$$
which implies $\|I_{\lambda,\beta}'(u)\|\leq \|u-T_{\lambda,\beta}(u)\|(1+C_1\|u\|^2+C_2\|u\|^{2\alpha})$ for all $u\in E$.
\qed
\begin{lemma}\label{Lem:T3}
For fixed $(\lambda,\beta)\in(0,1]\times(0,1]$ and for $c<d$ and $\tau>0$,
there exists $\delta>0$ (which depends on $\lambda$ and $\beta$) such that $\|u-T_{\lambda,\beta}(u)\|\geq \delta$
if $u\in E$, $I_{\lambda,\beta}(u)\in [c,d]$
and $\|I_{\lambda,\beta}'(u)\|\geq\tau$.
\end{lemma}
\Proof
Fix $\gamma\in(4,r)$, then for $u\in E$, we have
\begin{equation*}\label{eqn:T0}
\aligned
&I_{\lambda,\beta}(u)-\frac{1}{\gamma}\langle u,u-T_{\lambda,\beta}(u)\rangle\\
&=\frac{\gamma-2}{2\gamma}\|u\|^2
+\frac{b}{\gamma}\int_{\R^3}|\nabla u|^2\int_{\R^3}(\nabla u-\nabla T_{\lambda,\beta}(u))\nabla u\\
&+\lambda\frac{\gamma-2(1+\alpha)}{2\gamma(1+\alpha)}\|u\|_2^{2\alpha+2}+\frac{\lambda}{\gamma}\|u\|_2^{2\alpha}\int_{\R^3}u(u-T_{\lambda,\beta}(u))\\
&+\int_{\R^3}(\frac{1}{\gamma}f(u)u-F(u))+\frac{\gamma-4}{4\gamma}b\left(\int_{\R^3}|\nabla u|^2\right)^2
+\frac{(r-\gamma)\beta}{r\gamma}\int_{\R^3}|u|^r.
\endaligned
\end{equation*}
It follows from (\ref{eqn:f}) that for any $\epsilon>0$, there exists $C_\epsilon>0$ such that
\begin{equation*}\label{eqn:T1}
\aligned
&|I_{\lambda,\beta}(u)|+\frac{1}{\gamma}\|u\|\|u-T_{\lambda,\beta}(u)\|\\
&\geq(\frac{\gamma-2}{2\gamma}-\epsilon C)\|u\|^2
+\frac{b}{\gamma}\int_{\R^3}|\nabla u|^2\int_{\R^3}(\nabla u-\nabla T_{\lambda,\beta}(u))\nabla u\\
&+\frac{\gamma-4}{4\gamma}b\left(\int_{\R^3}|\nabla u|^2\right)^2+\frac{r-\gamma}{r\gamma}\beta\int_{\R^3}|u|^r
+\lambda\frac{\gamma-2(1+\alpha)}{2\gamma(1+\alpha)}\|u\|_2^{2\alpha+2}\\
&-C_\epsilon\|u\|_p^p+\frac{\lambda}{\gamma}\|u\|_2^{2\alpha}\int_{\R^3}u(u-T_{\lambda,\beta}(u)).
\endaligned
\end{equation*}
Then,
\begin{equation}\label{eqn:T2}
\aligned
&\|u\|^2+b\left(\int_{\R^3}|\nabla u|^2\right)^2+\beta\|u\|_r^r+\lambda\|u\|_2^{2\alpha+2}-C_\epsilon\|u\|_p^p\\
&\le C(|I_{\lambda,\beta}(u)|+\|u\|\|u- T_{\lambda,\beta}(u)\|
+\frac{\lambda}{\gamma}\|u\|_2^{2\alpha}\int_{\R^3}|u||u-T_{\lambda,\beta}(u)|\\
&+\frac{b}{\gamma}\int_{\R^3}|\nabla u|^2
\int_{\R^3}|\nabla u-\nabla T_{\lambda,\beta}(u)||\nabla u|).
\endaligned
\end{equation}
By H\"{o}lder's inequality and Sobolev's inequality, we have
\begin{equation}\label{eqn:T3}
\aligned
&\frac{b}{\gamma}\int_{\R^3}|\nabla u|^2\int_{\R^3}|\nabla u-\nabla T_{\lambda,\beta}(u)||\nabla u|
\leq C\int_{\R^3}|\nabla u|^2\|u\|\|u- T_{\lambda,\beta}(u)\|,\\
&\text{ and } \frac{\lambda}{\gamma}\|u\|_2^{2\alpha}\int_{\R^3}u(u-T_{\lambda,\beta}(u))\leq C\|u\|_2^{2\alpha}\|u\|\|u- T_{\lambda,\beta}(u)\|.
\endaligned
\end{equation}
By (\ref{eqn:T2}) and (\ref{eqn:T3}) and Young's inequality, we get
\begin{equation}\label{eqn:T4}
\aligned
&\|u\|^2+b\left(\int_{\R^3}|\nabla u|^2\right)^2+\beta\|u\|_r^r+\lambda\|u\|_2^{2\alpha+2}-C_\epsilon\|u\|_p^p\\
&\le C(|I_{\lambda,\beta}(u)|+\|u\|\|u- T_{\lambda,\beta}(u)\|+\|u\|^2\|u- T_{\lambda,\beta}(u)\|^2+\|u\|_2^{4\alpha}).
\endaligned
\end{equation}
Assume on the contrary that there exists $\{u_n\}\subset E$ with $I_{\lambda,\beta}(u_n)\in[c,d]$ and
$\|I'_{\lambda,\beta}(u_n)\|\geq \tau$ such that
$$
\|u_n-T_{\lambda,\beta}(u_n)\|\rightarrow0\quad\text{as}\,\,n\rightarrow\infty,
$$
then it follows from (\ref{eqn:T4}) that
\begin{equation}\label{eqn:T5}
\|u_n\|^2+b\left(\int_{\R^3}|\nabla u_n|^2\right)^2+\beta\|u_n\|_r^r+\lambda\|u_n\|_2^{2\alpha+2}-C_\epsilon\|u_n\|_p^p\leq C(1+\|u_n\|_2^{4\alpha})
\end{equation}
for large $n$. Now we claim that $\{u_n\}$ is a bounded sequence in $E$.
Otherwise, assume $\|u_n\|\rightarrow\infty$, then by (\ref{eqn:T5}) we have
\begin{equation}\label{eqn:T6}
\|u_n\|^2+b\left(\int_{\R^3}|\nabla u_n|^2\right)^2+\beta\|u_n\|_r^r+\lambda\|u_n\|_2^{2\alpha+2}-C_\epsilon\|u_n\|_p^p\leq C.
\end{equation}
Note that, for any $A_1>0$, we can choose $A_2>0$ such that
$$
t^{1+\alpha}>A_1t-A_2.
$$
Applying this with $t=\|u_n\|_2^2$, then by (\ref{eqn:T6}) we have
\begin{equation}\label{eqn:T7}
\|u_n\|^2+b\left(\int_{\R^3}|\nabla u_n|^2\right)^2+\int_{\R^3}(\beta|u_n|^r+\lambda A_1|u_n|^{2}-C_\epsilon|u_n|^p)-A_2\leq C.
\end{equation}
Since $2<p<r$, we take $A_1$ large enough such that the function $\lambda A_1|t|^{2}+\beta|t|^r-C_\epsilon |t|^p>0$ for any $t\in\R$.
Then (\ref{eqn:T7}) implies a contradiction. So our claim is true, that is, $\{u_n\}$ is a bounded sequence
in $E$ for any fixed $(\lambda,\beta)\in(0,1]\times(0,1]$.
The claim combined with Lemma \ref{Lem:T2} implies $\|I'_{\lambda,\beta}(u_n)\|\rightarrow0$ as $n\rightarrow\infty$, which is a contradiction.
The proof is complete.
\qed

\subsection{Invariant subsets of descending flow}

In order to obtain sign-changing solutions, we define the positive and negative cones by
$$
P^+:=\{u\in E:\, u\geq0\}\quad\text{and}\quad P^-:=\{u\in E:\,u\leq0\},
$$
respectively. For $\epsilon>0$, set
$$
P_\epsilon^+:=\{u\in E:\, dist(u,P^+)<\epsilon\}\quad\text{and}\quad P_\epsilon^-:=\{u\in E:\,dist(u,P^-)<\epsilon\},
$$
where $dist(u,P^\pm)=\inf\limits_{v\in P^\pm}\|u-v\|$. Clearly, $P_\epsilon^-=-P_\epsilon^+$. Let $W=P_\epsilon^+\cup P_\epsilon^-$.
It is easy to check that $W$ is an open and symmetric subset of $E$ and $E\setminus {W}$ contains only sign-changing functions.

 We denote by $K$ the set of critical points of $I_{\lambda,\beta}$, that is,
$K=\{u\in E:\,\,I_{\lambda,\beta}'(u)=0\}$ and $E_0:=E\setminus{K}$.
For $c\in\R$, define $K_c=\{u\in E:I_{\lambda,\beta}(u)=c,I_{\lambda,\beta}'(u)=0\}$ and
$I_{\lambda,\beta}^c=\{u\in E:I_{\lambda,\beta}(u)\leq c\}$.

In the following, we will show that, for $\epsilon$ small enough,
all sign-changing solutions to (K$_{\lambda,\beta}$) are contained in $E\setminus {W}$.

\begin{lemma}\label{Lem:C1}
 There exists $\epsilon_0>0$ such that for $\epsilon\in(0,\epsilon_0)$,
\item[\rm (1) ] $T_{\lambda,\beta}(\partial P_\epsilon^-)\subset P_\epsilon^-$ and every nontrivial solution $u\in P_\epsilon^-$ is negative,
\item[\rm (2) ] $T_{\lambda,\beta}(\partial P_\epsilon^+)\subset P_\epsilon^+$ and every nontrivial solution $u\in P_\epsilon^+$ is positive.
\end{lemma}
\Proof
We only prove $T_{\lambda,\beta}(\partial P_\epsilon^-)\subset P_\epsilon^-$, and the other case is similar.
For $u\in E$, define $v:=T_{\lambda,\beta}(u)$. Since $\text{dist}(v,P^-)\leq \|v^+\|$, by Sobolev's inequality
and ($f_1$)-($f_2$),  for any $\varepsilon>0$, there exists $C_\varepsilon>0$ such that
$$
\aligned
&\text{dist}(v,P^-)\|v^+\|\leq\|v^+\|^2=\langle v,v^+\rangle\\
&\leq\int_{\R^3}f(u)v^+-b\int_{\R^3}|\nabla u|^2\int_{\R^3}\nabla v\nabla v^++\int_{\R^3}|u|^{r-2}uv^+
-\lambda\|u\|_2^{2\alpha}\int_{\R^3}vv^+\\
&\leq\int_{\R^3}f(u)v^++\int_{\R^3}|u|^{r-2}uv^+\\
&\leq\int_{\R^3}f(u^+)v^++\int_{\R^3}|u^+|^{r-2}u^+v^+\\
&\leq\int_{\R^3}(\varepsilon u^+v^++C_\varepsilon|u^+|^{p-1}v^+)+\int_{\R^3}|u^+|^{r-2}u^+v^+\\
&\leq C[\varepsilon\text{dist}(u,P^-)+C_\varepsilon\text{dist}(u,P^-)^{p-1}+\text{dist}(u,P^-)^{r-1}]\|v^+\|,
\endaligned
$$
which further implies that
$$
\text{dist}(v,P^-)\leq C[\varepsilon\text{dist}(u,P^-)+C_\varepsilon\text{dist}(u,P^-)^{p-1}+\text{dist}(u,P^-)^{r-1}].
$$
There exists $\epsilon_0>0$ such that for $\epsilon\in(0,\epsilon_0)$,
$$
\text{dist}(T_{\lambda,\beta}(u),P^-)=\text{dist}(v,P^-)<\epsilon.
$$
Therefore, we have $T_{\lambda,\beta}(u)\in P_\epsilon^-$ for any $u\in P_\epsilon^-$.
\qed

Since the operator $T_{\lambda,\beta}$ is not locally Lipschitz continuous,
we need to construct a locally Lipschitz continuous vector field which inherits its properties.
Arguing as the proof of Lemma 2.1 in \cite{Bartsch05}, we have
\begin{lemma}\label{Lem:local}
There exists a locally Lipschitz continuous operator $B_{\lambda,\beta}:\,E\rightarrow E$ such that
\begin{itemize}
\item[(i)] $\langle T'_{\lambda,\beta}(u),u-B_{\lambda,\beta}\rangle\geq\frac{1}{2}\|u-T_{\lambda,\beta}(u)\|^2$;
\item[ (ii)]$\frac{1}{2}\|u-B_{\lambda,\beta}(u)\|^2\leq\|u-T_{\lambda,\beta}(u)\|^2\leq2\|u-B_{\lambda,\beta}(u)\|^2$;
\item[ (iii)]$T_{\lambda,\beta}(\partial P_\epsilon^\pm)\subset P_\epsilon^\pm$, $\forall \epsilon\in(0,\epsilon_0)$;
\item[ (iv)]if $I_{\lambda,\beta}$ is even, then $B_{\lambda,\beta}$ is odd.
\end{itemize}
\end{lemma}
In what follows, we verifies that the functional $I_{\lambda,\beta}$ satisfies (PS)-condition.

\begin{lemma}\label{Lem:PS}
Assume that there exist $\{u_n\}\subset E$ and $c\in\R$ such that $I_{\lambda,\beta}(u_n)\rightarrow c$ and
$I'_{\lambda,\beta}(u_n)\rightarrow 0$ for any fixed $(\lambda,\beta)\in(0,1]\times(0,1]$ as $n\rightarrow\infty$, then there exists a convergence subsequence of $\{u_n\}$,
 denoted by $\{u_n\}$ for simplicity, such that $u_n\rightarrow u$ in $E$ for some $u\in E$.
\end{lemma}
\Proof
For $\gamma\in(4,p)$, we have
$$
\aligned
&\gamma I_{\lambda,\beta}(u_n)-\langle I'_{\lambda,\beta}(u_n),u_n\rangle\\
&=\frac{\gamma-2}{2}\|u_n\|^2+\frac{b(\gamma-4)}{4}\left(\int_{\R^3}|\nabla u_n|^2\right)^2\\
&+\lambda\frac{\gamma-2(1+\alpha)}{2(1+\alpha)}\|u_n\|_2^{2(1+\alpha)}+\int_{\R^3}( f(u_n)u_n-\gamma F(u_n))+\beta\frac{r-\gamma}{r}\int_{\R^3}|u_n|^r.
\endaligned
$$
It follows from ($f_1$)-($f_2$) and Sobolev's inequality that for any $\varepsilon>0$, there exists $C_\varepsilon>0$ such that
$$
\aligned
&\gamma |I_{\lambda,\beta}(u_n)|+\|I'_{\lambda,\beta}(u_n)\|\|u_n\|\\
&=(\frac{\gamma-2}{2}-\varepsilon C)\|u_n\|^2+\frac{b(\gamma-4)}{4}\left(\int_{\R^3}|\nabla u_n|^2\right)^2-C_\varepsilon\|u_n\|_p^p\\
&+\lambda\frac{\gamma-2(1+\alpha)}{2(1+\alpha)}\|u_n\|_2^{2(1+\alpha)}+\beta\frac{r-\gamma}{r}\int_{\R^3}|u_n|^r.
\endaligned
$$
Arguing as the proof of Lemma \ref{Lem:T3}, one sees that $\{u_n\}$ is bounded in $E$. Up to subsequence, we assume that there exists $u\in E$ such that
$$
\aligned
&u_n\rightharpoonup u\,\,\text{weakly\,in}\,\,E, \text{ and }\\
& u_n\rightarrow u\,\,\text{strongly\,in}\,\,L^q(\R^3)\,\text{for}\,q\in(2,6).
\endaligned
$$
Note that
\begin{equation}\label{eqn:PS1}
\aligned
&\langle I'_{\lambda,\beta}(u_n)-I'_{\lambda,\beta}(u),u_n-u\rangle\\
&=\|u_n-u\|^2+b\int_{\R^3}|\nabla u_n|^2\int_{\R^3}|\nabla(u_n-u)|^2\\
& +b(\int_{\R^3}|\nabla u_n|^2-\int_{\R^3}|\nabla u|^2)\int_{\R^3}\nabla u\nabla(u_n-u)-\int_{\R^3}(f(u_n)-f(u))(u_n-u)\\
&+\lambda\|u_n\|_2^{2\alpha}\int_{\R^3}(u_n-u)^2+\lambda(\|u_n\|_2^{2\alpha}-\|u\|_2^{2\alpha})\int_{\R^3}u(u_n-u)\\
&-\beta\int_{\R^3}(|u_n|^{r-2}u_n-|u|^{r-2}u)(u_n-u).
\endaligned
\end{equation}
By the boundedness of $\{u_n\}$ in $E$, one has
$$
\aligned
&b(\int_{\R^3}|\nabla u_n|^2-\int_{\R^3}|\nabla u|^2)\int_{\R^3}\nabla u\nabla(u_n-u)\rightarrow0,\\
&\lambda(\|u_n\|_2^{2\alpha}-\|u\|_2^{2\alpha})\int_{\R^3}u(u_n-u)\rightarrow0,
\quad \text{as}\,\,n\rightarrow\infty.
\endaligned
$$
 Moreover, for any $\varepsilon>0$, one has
$$
\aligned
\int_{\R^3}(f(u_n)-f(u))(u_n-u)&\leq\int_{\R^3}[\varepsilon(|u_n|+|u|)+C_\varepsilon(|u_n|^{p-}+|u_n|^{p-1})]|u_n-u|\\
&\leq\varepsilon C+C_\varepsilon(\|u_n\|_p^{p-1}+\|u\|_p^{p-1})\|u_n-u\|_p\rightarrow0
\endaligned
$$
as $n\rightarrow\infty$. Similarly, we also have
$$
\beta\int_{\R^3}(|u_n|^{r-2}u_n-|u|^{r-2}u)(u_n-u)\rightarrow0,\quad \text{as}\,\,n\rightarrow\infty.
$$
Based on the above facts, from (\ref{eqn:PS1}) we deduce that $u_n\rightarrow u$ in $E$ as $n\rightarrow\infty$.
\qed

Here we give a deformation lemma to functional $I_{\lambda,\beta}$
whose proof is almost the same as that of Lemma 3.6 in \cite{Liuz16}.
\begin{lemma}\label{Lem:XB}(Deformation lemma)
Let $S\subset E$ and $c\in\R$ such that
$$
\forall u\in I^{-1}_{\lambda,\beta}([c-2\epsilon_0,c+2\epsilon_0])\cap S_{2\delta},\,\,\,\|I'_{\lambda,\beta}(u)\|\geq \epsilon_0,
$$
where $\epsilon_0$ was given in Lemma \ref{Lem:C1} and $S_{2\delta}:=\{u\in S:\,\,\text{dist}(u,S)<2\delta\}$. Then for
$\epsilon_1\in(0,\epsilon_0)$ there exists $\eta\in C([0,1]\times E,E)$ such that
\begin{itemize}
\item[(i)] $\eta(t,u)=u$ if $t=0$ or if $u\not\in I^{-1}_{\lambda,\beta}([c-2\epsilon_1,c+2\epsilon_1])$;
\item[ (ii)]$\eta(1,I^{c+\epsilon_1}_{\lambda,\beta}\cap S)\subset I^{c-\epsilon_1}_{\lambda,\beta}$;
\item[ (iii)] $I_{\lambda,\beta}(\eta(\cdot,u))$ is not increasing for all $u\in E$;
\item[ (iv)]$\eta(t, \overline{P_\epsilon^+})\subset \overline{P_\epsilon^+}$, $\eta(t,\overline{P_\epsilon^-})\subset \overline{P_\epsilon^-}$, $\forall t\in [0,1]$;
\item[ (v)] if $f$ is odd, then $\eta(t,\cdot)$ is odd $\forall t\in[0,1]$.
\end{itemize}
\end{lemma}

Now we introduce a critical point theorem (see \cite{Liuj15}). For more details, we let
$P,Q\subset E$ be open sets, $M=P\cap Q$,
$\Sigma=\partial P\cap\partial Q$ and $W=P\cup Q$.

\begin{definition} {\rm (see \cite{Liuj15})} \label{Def:1}
$\{P,Q\}$ is called an admissible family of invariant sets
with respect to $J$ at level $c$, provided that the following deformation property
holds: if $K_c\setminus{W }=\emptyset$, then, there exists $\epsilon_0 > 0$ such that for any
$\epsilon\in (0,\epsilon_0)$, there
exists $\eta\in C(E,E)$ satisfying
\begin{itemize}
\item[\rm (1) ]
$\eta(\bar{P})\subset\bar{P}$, $\eta(\bar{Q})\subset\bar{Q}$;
\item[\rm (2) ]$\eta|_{J^{c-2\epsilon}}=id$;
\item[\rm (3) ] $\eta(J^{c+\epsilon}\setminus{W})\subset J^{c-\epsilon}$.
\end{itemize}
\end{definition}

\begin{theorem}{\rm (see \cite{Liuj15})}\label{Thm:xjl}
 Assume that $\{P,Q\}$ is an admissible family of invariant
sets with respect to $J$ at any level $c\geq c_* := \inf_{u\in\Sigma}J(u)$ and there exists a map
$\psi_0 :\triangle\rightarrow E$ satisfying
\begin{itemize}
\item[\rm (1) ]
$\psi_0(\partial_1\triangle)\subset P$ and $\psi_0(\partial_2\triangle)\subset Q$,
\item[\rm (2) ]$\psi_0(\partial_0\triangle)\cap M=\emptyset$,
\item[\rm (3) ] $\sup\limits_{u\in\psi_0(\partial_0\triangle)}J(u)<c_*$,
\end{itemize}
where $\triangle=\{(t_1,t_2)\in\R^2:t_1,t_2>0,\,t_1+t_2\leq1\}$, $\partial_1\triangle=\{0\}\times[0,1]$,
$\partial_2\triangle=[0,1]\times\{0\}$ and  $\partial_0\triangle=\{(t_1,t_2)\in\R^2:t_1,t_2\geq0,\,t_1+t_2=1\}$.
Define
$$
c=\inf\limits_{\psi\in\Gamma}\sup\limits_{u\in\psi(\triangle)\setminus{W}}J(u),
$$
where $\Gamma:=\{\psi\in C(\triangle,E):\,\psi(\partial_1\triangle)\subset P,\,
\psi(\partial_2\triangle)\subset Q,\,\psi|_{\partial_0\triangle}=\psi_0|_{\partial_0\triangle}\}$.
Then $c\geq c_*$ and $K_c\setminus{W}\not=\emptyset$.
\end{theorem}

\subsection{Proof of Theorem \ref{Thm:existence}}
In order to employ Theorem \ref{Thm:xjl} to prove the existence of sign-changing solutions
to problem (K$_{\lambda,\beta}$), we take
$P=P_\epsilon^+,Q=P_\epsilon^-$ and $J=I_{\lambda,\beta}$. We need to prove that $\{P_\epsilon^+,P_\epsilon^-\}$
is an admissible family of invariant sets for the
functional $I_{\lambda,\beta}$ at any level $c\in \R$.
Moreover, $K_c\subset W$ if $K_c \setminus{W}=\emptyset$. Since
the functional $I_{\lambda,\beta}$ satisfies the (PS)-condition, $K_c$ is compact. Thus, $2\delta:= dist(K_c; \partial W)>0$.

\begin{lemma}\label{Lem:qfanshu}
For $q\in[2,6]$, there exits $m>0$ independent of $\epsilon$ such that
$\|u\|_q\leq m\epsilon$ for $u\in M=P_\epsilon^+\cap P_\epsilon^-$.
\end{lemma}
\Proof
For any fixed $u\in M$, we have
$$
\|u^\pm\|_q=\inf\limits_{v\in P^{\mp}}\|u-v\|_q\leq C\inf\limits_{v\in P^\mp}\|u-v\|\leq C\text{dist}(u,P^{\mp}).
$$
Then $\|u\|_q\leq m\epsilon$ for $u\in M$.
\qed

\begin{lemma}\label{Lem:xiaozhi}
If $\epsilon>0$ is small enough, then $I_{\lambda,\beta}(u)\geq\frac{\epsilon^2}{4}$ for all $u\in \Sigma=\partial P_\epsilon^+\cap\partial P_\epsilon^-$,
that is, $c_*\geq\frac{\epsilon^2}{4}$.
\end{lemma}
\Proof
For any fixed $u\in \partial P_\epsilon^+\cap\partial P_\epsilon^-$, we have
$\|u^+\|\geq\text{dist}(u,P^{-})=\epsilon$.
By Lemma \ref{Lem:qfanshu} and ($f_1$)-($f_2$),  for any $\varepsilon>0$, there exists $C_\varepsilon>0$ such that
$$
\aligned
I_{\lambda,\beta}(u)&=\frac{1}{2}\|u\|^2+\frac{b}{4}\left(\int_{\R^3}|\nabla u|^2\right)^2+\frac{\lambda}{2(1+\alpha)}\|u\|_2^{2(1+\alpha)}-\int_{\R^3}F(u)-\frac{\beta}{r}\int_{\R^3}|u|^r\\
&\geq\frac{1}{2}\epsilon^2-\frac{\varepsilon}{2}\|u\|_2^2-\frac{C_\varepsilon}{p}\|u\|_p^p-\frac{\beta}{r}\|u\|_r^r\\
&\geq\frac{1}{2}\epsilon^2-\frac{\varepsilon}{2}\epsilon^2-\frac{C_\varepsilon}{p}\epsilon^p-
\frac{\beta}{r}\epsilon^r\\
&\geq\frac{\epsilon^2}{4}
\endaligned
$$
for $\epsilon$ and $\varepsilon$ small enough.
\qed

\textbf{Proof of Theorem \ref{Thm:existence}.} We use Theorem \ref{Thm:xjl} to prove
the existence of sign-changing
solutions to problem (K$_{\lambda,\beta}$). Let $X=E$, $P=P_\epsilon^+$, $Q=P_\epsilon^-$ and $J=I_{\lambda,\beta}$.
Take $S=E\setminus{W}$ in Lemma \ref{Lem:XB}, then we can easily deduce that $\{P_\epsilon^+,P_\epsilon^-\}$ is an admissible
family of invariant sets for the functional $I_{\lambda,\beta}$ at any level $c\in\R$.\\
In what follows, we divide three steps to complete the proof.\\
\textbf{Step 1}. Choose $v_1,v_2\in C_0^{\infty}(B_1(0))$ such that $\text{supp}(v_1)\cap\text{supp}(v_2)=\emptyset$ and
$v_1<0,v_2>0$, where $B_r(0):=\{x\in\R^3:\,|x|<r\}$. For $(t,s)\in \triangle$, let
$$
\varphi_0(t,s):=R^2[tv_1(R\cdot)+sv_2(R\cdot)],
$$
where $R>0$ will be determined later. Obviously, for $t,s\in[0,1]$,
$\varphi_0(0,s)(\cdot)=Rsv_2(R^2\cdot)\in P_\epsilon^+$
and $\varphi_0(t,0)(\cdot)=R^2tv_1(R\cdot)\in P_\epsilon^-$.
It follows from Lemma \ref{Lem:xiaozhi} that, for small $\epsilon>0$,
$$
I_{\lambda,\beta}(u)\geq \frac{\epsilon^2}{4}\quad\text{for\,all}\,u\in \Sigma=\partial P_\epsilon^+\cap\partial P_\epsilon^-,\, (\lambda,\beta)\in(0,1]\times(0,1].
$$
Hence $c_*=\inf_{u\in\Sigma}I_{\lambda,\beta}(u)\geq\frac{\epsilon^2}{4}$ for any $(\lambda,\beta)\in(0,1]\times(0,1]$.
Let $u_t=\varphi_0(t,1-t)$ for $t\in [0,1]$. Observe that
$$
\rho=\min\{\|tv_1+(1-t)v_2\|_2:\,\,0\leq t\leq 1\}>0,
$$
then $\|u_t\|_2^2\geq\rho R$ for $u\in\varphi_0(\partial_0\triangle)$. It follows from Lemma \ref{Lem:qfanshu} that
$\varphi_0(\partial_0\triangle)\cap P_\epsilon^+\cap P_\epsilon^-=\emptyset$.
 A direct computation shows that
\begin{equation}\label{eqn:PS1--}
\aligned
&\int_{\R^3}|\nabla u_t|^2=R^3\int_{\R^3}(t^2|\nabla v_1|^2+(1-t)^2|\nabla v_2|^2)=:R^3B(t),\\
&\int_{\R^3}V(x)|u_t|^2\leq R \max_{x\in B_1(0)}V(x)\int_{\R^3}(t^q|v_1|^q+(1-t)^q|v_2|^q)=:RB_2(t),\\
&\int_{\R^3}|u_t|^q=R^{2q-3}\int_{\R^3}(t^q|v_1|^q+(1-t)^q|v_2|^q)=:R^{2q-3}B_q(t)\,\,\text{for}\,\,q\in(2,6],\\
&\left(\int_{\R^3}|u_t|^2\right)^{1+\alpha}=R^{(1+\alpha)}\left(\int_{\R^3}(t^2|v_1|^2+(1-t)^2|v_2|^2)\right)^{1+\alpha}=:R^{(1+\alpha)}\bar{B}(t).
\endaligned
\end{equation}
Since $F(t)\geq C_3|t|^\mu-C_4$ for any $t\in\R$
and some positive constants $C_3,C_4$, we have
$$
\aligned
I_{\lambda,\beta}(u_t)=&\frac{a}{2}\int_{\R^3}|\nabla u_t|^2+\frac{1}{2}\int_{\R^3}V(x)|u_t|^2
+\frac{1}{2(1+\alpha)}\|u_t\|_2^{2(1+\alpha)}\\
&+\frac{b}{4}\left(\int_{\R^3}|\nabla u_t|^2\right)^2-\int_{B_{R^{-1}}(0)}F(u_t)-\frac{\beta}{r}\int_{\R^3}|u_t|^r\\
<&\frac{aR^3}{2}B(t)+\frac{R}{2}B_2(t)+\frac{bR^{6}}{4}B^2(t)+\frac{R^{(1+\alpha)}}{2(1+\alpha)}\bar{B}(t)\\
&-C_3R^{2\mu-3}B_\mu(t)+CC_4R^{-3}-\frac{R^{2r-3}}{r}B_r(t).
\endaligned
$$
Since $r\in(\max\{p,\frac{9}{2}\},6)$, one sees that $I_{\lambda,\beta}(u_t)\rightarrow-\infty$
as $R\rightarrow+\infty$ for any fixed
$(\alpha,\beta)\in(0,1]\times(0,1]$. Hence we can choose $R$ large enough such that
$$
\sup\limits_{u\in\varphi_0(\partial_0\triangle)}I_{\lambda,\beta}(u)<c_*:=\inf\limits_{u\in\Sigma}I_{\lambda,\beta}(u).
$$
Since $I_{\lambda,\beta}$ satisfies the assumptions of Theorem \ref{Thm:xjl}, the number
$$
c_{\lambda,\beta}=\inf\limits_{\varphi\in\Gamma}\sup\limits_{u\in\varphi(\triangle)\setminus{W}}I_{\lambda,\beta}(u)
$$
is a critical value of $I_{\lambda,\beta}$ satisfying $c_{\lambda,\beta}\geq c_*$. Therefore, there exists $u_{\lambda,\beta}\in E\setminus{(P_\epsilon^+\cup P_\epsilon^-)}$ such that
$I_{\lambda,\beta}(u_{\lambda,\beta})=c_{\lambda,\beta}$ and $I_{\lambda,\beta}'(u_{\lambda,\beta})=0$ for $(\lambda,\beta)\in(0,1]\times(0,1]$.

\textbf{Step 2.}
Passing to the limit as $\lambda\rightarrow0$ and $\beta\rightarrow0$. According to the definition of $c_{\lambda,\beta}$,
we see that for any $(\lambda,\beta)\in(0,1]\times(0,1]$,
\begin{equation}\label{eqn:fun0}
c_{\lambda,\beta}\leq C_R:=\sup\limits_{u\in\varphi_0(\triangle)}I_{1,0}(u)<\infty,
\end{equation}
where $C_R$ is independent of $(\lambda,\beta)\in(0,1]\times(0,1]$.
Without loss of generality, we set $\lambda=\beta$. Choosing a sequence $\{\lambda_n\}\subset(0,1]$ satisfying $\lambda_n\rightarrow0^+$, then
we find a sequence of sign-changing critical points $\{u_{\lambda_n}\}$ (still denoted by $\{u_n\}$ for simplicity) of $I_{\lambda_n,\beta_n}$ and
$I_{\lambda_n,\beta_n}(u_n)=c_{\lambda_n,\beta_n}$.
Now we show that $\{u_n\}$ is bounded in $E$. By the definition of $I_{\lambda,\beta}$, we have
\begin{equation}\label{eqn:fun1}
\aligned
c_{\lambda_n,\beta_n}=\frac{a}{2}\int_{\R^3}|\nabla u_n|^2&+\frac{1}{2}\int_{\R^3}V(x)u_n^2
+\frac{\lambda}{2(1+\alpha)}\|u_n\|_2^{2(1+\alpha)}\\
&+\frac{b}{4}\left(\int_{\R^3}|\nabla u_n|^2\right)^2-\int_{\R^3}F(u_n)-\frac{\beta}{r}\int_{\R^3}|u_n|^r
\endaligned
\end{equation}
and
\begin{equation}\label{eqn:fun2}
\aligned
0=a\int_{\R^3}|\nabla u_n|^2&+\int_{\R^3}V(x)u_n^2
+\lambda\|u_n\|_2^{2(1+\alpha)}\\
&+b\left(\int_{\R^3}|\nabla u_n|^2\right)^2-\int_{\R^3}f(u_n)u_n-\beta\int_{\R^3}|u_n|^r.
\endaligned
\end{equation}
Moreover, from Lemma \ref{Lem:pohozaev},  the following identity holds
\begin{equation}\label{eqn:fun3}
\aligned
& \frac{a}{2}\int_{\R^3}|\nabla u_n|^2+\frac{3}{2}\int_{\R^3}V(x)u_n^2+\frac{1}{2}\int_{\R^3}(\nabla
V(x),x)u_n^2\\+&\frac{b}{2}\left(\int_{\R^3}|\nabla u_n|^2\right)^2
+\frac{3\lambda}{2}\left(\int_{\R^3}u_n^2\right)^{1+\alpha}-3\int_{\R^3}(F(u_n)+\frac{\beta}{r}|u_n|^r)=0.
\endaligned
\end{equation}
Multiplying (\ref{eqn:fun1}), (\ref{eqn:fun2}) and (\ref{eqn:fun3}) by $4$, $-\frac{1}{\mu}$ and  $-1$ respectively and adding them up,
we get
\begin{equation*}\label{eqn:fun4}
\aligned
4c_{\lambda_n,\beta_n}=&a\frac{3\mu-2}{2\mu}\int_{\R^3}|\nabla u_n|^2
+\frac{\mu-2}{2\mu}\int_{\R^3}V(x)u_n^2
-\frac{1}{2}\int_{\R^3}(\nabla
V(x),x)u_n^2\\
&+\frac{\mu-2}{2\mu}b\left(\int_{\R^3}|\nabla u_n|^2\right)^2+\lambda\frac{\mu-2-3\mu\alpha-2\alpha}{2\mu(1+\alpha)}\|u_n\|_2^{2(1+\alpha)}\\
&+\int_{\R^3}(\frac{1}{\mu}f(u_n)u_n-F(u_n))+\beta\frac{r-\mu}{\mu r}\int_{\R^3}|u_n|^r.
\endaligned
\end{equation*}
Since $\alpha<\frac{\mu-2}{3\mu+2}$ and $\mu>2$, it follows from ($V_2$), ($f_3$) and (\ref{eqn:fun0}) that
\begin{equation*}\label{eqn:fun5}
4C_R>a\frac{3\mu-2}{2\mu}\int_{\R^3}|\nabla u_n|^2+\frac{\mu-2}{2\mu}b\left(\int_{\R^3}|\nabla u_n|^2\right)^2,
\end{equation*}
which implies that there exists $C_5>0$ independent of $\lambda,\beta$ such that
\begin{equation}\label{eqn:fun6}
\int_{\R^3}|\nabla u_n|^2<C_5.
\end{equation}
Moreover, combining (\ref{eqn:fun0}), (\ref{eqn:fun1}) and hypotheses (V$_1$), ($f_1$) and ($f_2$), we obtain that
for small $\varepsilon>0$, there exists $C_\varepsilon>0$ such that
\begin{equation}\label{eqn:fun7}
\aligned
C_R&>\frac{a}{2}\int_{\R^3}|\nabla u_n|^2+\frac{1}{2}\int_{\R^3}V(x)u_n^2
-\int_{\R^3}F(u_n)-\frac{\beta}{r}\int_{\R^3}|u_n|^r\\
&>\frac{1-\varepsilon}{2}\int_{\R^3}V(x)u_n^2-C_\varepsilon\int_{\R^3}u_n^6-\frac{1}{r}\int_{\R^3}|u_n|^r\\
&>\frac{1-\varepsilon}{2}\int_{\R^3}V(x)u_n^2-C_\varepsilon S^{-3}\left(\int_{\R^3}|\nabla u_n|^2\right)^3-\frac{1}{r}\int_{\R^3}|u_n|^r.
\endaligned
\end{equation}
From interpolation inequality, Sobolev's inequality and Young's inequality, we deduce that for $\varepsilon>0$,
there exists $\bar{C}_\varepsilon>0$
such that
\begin{equation}\label{eqn:fun8}
\aligned
\int_{\R^3}|u_n|^r&\leq \left(\int_{\R^3}u_n^2\right)^{\frac{6-r}{4}}\left(\int_{\R^3}|u_n|^6\right)^{\frac{r-2}{4}}\\
&\leq \varepsilon\left(\int_{\R^3}u_n^2\right)^{\frac{6-r}{2}}+\bar{C}_\varepsilon\left(\int_{\R^3}|u_n|^6\right)^{\frac{r-2}{2}}\\
&\leq \varepsilon\left(\int_{\R^3}u_n^2\right)^{\frac{6-r}{2}}+\bar{C}_\varepsilon S^{\frac{3(2-r)}{2}}
\left(\int_{\R^3}|\nabla u_n|^2\right)^{\frac{3(r-2)}{2}}.
\endaligned
\end{equation}
Combining (\ref{eqn:fun6}), (\ref{eqn:fun7}) and (\ref{eqn:fun8}), we immediately see that $\{u_n\}$ is bounded in $E$.
In view of (\ref{eqn:fun0}) and Lemma \ref{Lem:xiaozhi}, we have
$$
\aligned
\lim\limits_{n\rightarrow\infty}I(u_n)&=\lim\limits_{n\rightarrow\infty}\left(I_{\lambda_n,\beta_n}(u_n)-\frac{\lambda_n}{2(1+\alpha)}\|u_n\|_2^{2(1+\alpha)}
+\frac{\beta_n}{r}\int_{\R^3}|u_n|^r\right)\\
&=\lim\limits_{n\rightarrow\infty}c_{\lambda_n,\beta_n}= c^*>\frac{\epsilon^2}{4}.
\endaligned
$$
Moreover, for any $\psi\in C_0^\infty(\R^3)$,
$$
\lim\limits_{n\rightarrow\infty}I'(u_n)\psi=\lim\limits_{n\rightarrow\infty}
\left(I'_{\lambda_n,\beta_n}(u_n)\psi-\lambda_n\|u_n\|_2^{2\alpha}\int_{\R^3}u_n\psi
+\beta_n\int_{\R^3}|u_n|^{r-2}u_n\psi\right)=0.
$$
That is to say, $\{u_n\}$ is a bounded Palais-Smale sequence for $I$ at level $c^*$. Thus, there exists $u^*\in E$ such
that $u_n\rightharpoonup u^*$ weakly in $E$ and $u_n\rightarrow u^*$ strongly in $L^q(\R^3)$ for $q\in(2,6)$.
Similar argument of Lemma \ref{Lem:PS} lead to that $I'(u^*)=0$ and $u_n\rightarrow u^*$ strongly in $E$ as $n\rightarrow0$.
Thus, the fact that $u_n\in E\setminus(P_\epsilon^+\cup P_\epsilon^-)$ yields $u^*\in E\setminus(P_\epsilon^+\cup P_\epsilon^-)$ and then
$u^*$ is a sign-changing solution of (K).

\textbf{Step 3.} Define
$$
\bar{c}:=\inf\limits_{u\in\Theta}I(u),\quad\Theta:=\{u\in E\setminus\{0\},\,I'(u)=0,\,u^{\pm}\not\equiv0\}.
$$
Based on Step 2, we see that $\Theta\not=\emptyset$ and
$\bar{c}\leq c^*$, where $c^*$ is given in the Step 2.
 By the definition of $\bar{c}$, there exists
$\{u_n\}\subset E$ such that $I(u_n)\rightarrow \bar{c}$ and $I'(u_n)=0$. Using the earlier arguments,
we can prove that $\{u_n\}$ is bounded in $E$.
Arguing as in Lemma \ref{Lem:PS}, there exists a nontrivial
$u\in E$ such that $I(u)=\bar{c}$ and $I'(u)=0$.
Furthermore, we deduce from $\langle I'(u_n),u_n^\pm\rangle=0$ that for any $\varepsilon>0$,
there exists $C_\varepsilon>0$ such that
$$
\aligned
C(\|u_n^\pm\|_p^2+\int_{\R^3}|u^\pm_n|^2)\leq\|u_n^\pm\|^2&\leq\int_{\R^3}f(u_n)u_n^\pm=\int_{\R^3}f(u_n^\pm)u_n^\pm\\
&\leq \varepsilon\int_{\R^3}|u_n^\pm|^2+C_\varepsilon\int_{\R^3}|u_n^\pm|^p\\
&\leq \varepsilon \|u_n^\pm\|_2^2+C_\varepsilon \|u_n^\pm\|_p^p,
\endaligned
$$
which, together with the boundedness of $\{u_n\}$ in $E$, implies that $\|u_n^\pm\|_p\geq C$. Hence,
$\|u^\pm\|_p\geq C$, and then
$u$ is a ground state sign-changing solution of problem (K). The proof is complete.\qed

\section{Multiplicity}\setcounter{equation}{0}\label{sec4}
\renewcommand{\theequation}{4.\arabic{equation}}
In this section, we prove the existence of infinitely many sign-changing solutions to problem (K).

\subsection{Proof of Theorem \ref{Thm:many} (Multiplicity)}

In order to obtain infinitely many sign-changing solutions, we introduce an abstract critical point approach
developed by Liu et al \cite{Liuj15}, which we
recall below.
The notations from Section \ref{sec2} are still valid. Assume $G : E\rightarrow E$ is an
isometric involution, that is, $G^2=id$ and $d(Gx;Gy) = d(x;y)$ for $x$, $y\in E$. A subset $F\subset E$ is said to be symmetric if $Gx \in F$
for any $x \in F$. We
assume $J$ is $G$-invariant on $E$ in the sense that $J(Gx) = J(x)$ for any $x\in E$.
We also assume $Q = GP$. The genus of a closed symmetric subset $F$ of $E\setminus\{0\}$ is denoted
by $\gamma(F)$.
\begin{definition}{\rm (see \cite{Liuj15})}\label{Def:2}
$P$ is called a $G$-admissible invariant set with respect to
$J$ at level $c$, if the following deformation property holds: there exist $\epsilon_0>0$ and
a symmetric open neighborhood $N$ of $K_c\setminus{W}$ with $\gamma(\bar{N})< \infty$, such that for
$\epsilon\in(0,\epsilon_0)$ there exists $\eta\in C(E,E)$ satisfying
\begin{itemize}
\item[\rm (1) ]
$\eta(\bar{P})\subset\bar{P}$, $\eta(\bar{Q})\subset\bar{Q}$;
\item[\rm (2) ]$\eta\circ G=G\circ\eta$;
\item[\rm (3) ]$\eta|_{J^{c-2\epsilon}}=id$;
\item[\rm (4) ] $\eta(J^{c+\epsilon}\setminus{(N\cup W})\subset J^{c-\epsilon}$.
\end{itemize}
\end{definition}

\begin{theorem}{\rm (see \cite{Liuj15})}\label{Thm:xjl2}
 Assume that $P$ is a $G$-admissible invariant set with respect
to $J$ at any level $c\geq c_*:= \inf_{u\in\Sigma} J(u)$ and for any $n\in N$, there exists a
continuous map $\psi_n: B_n := \{x \in \R^n : |x|\leq1\}\rightarrow E$ satisfying
\begin{itemize}
\item[\rm (1) ] $\psi_n(0)\subset M:=P\cap Q$ and $\psi_n(-t)=G \psi_n(t)$ for $t\in B_n$;
\item[\rm (2) ]$\psi_n(\partial B_n)\cap M=\emptyset$;
\item[\rm (3) ] $\sup\limits_{u\in Fix_{G}\cup\psi_n(\partial B_n)}J(u)<c_*$, where $Fix_{G}:=\{u\in E;Gu=u\}$.
\end{itemize}
For $j\in N$, define
$$
c_j=\inf\limits_{B\in\Gamma_j}\sup\limits_{u\in B\setminus{W}}J(u),
$$
where
$$
\aligned
\Gamma_j:=\{&B\big{|}B=\psi(B_n\setminus{Y})for\, some\,\psi\in G_n, n\geq j,\\
&and\, open Y\subset B_n \,such\, that
-Y=Y\,and\, \gamma(\bar{Y})\leq n-j\}
\endaligned
$$
and
$$
G_n:=\{\psi|\psi\in C(B_n,E), \psi(-t)=G\psi(t)\,for\,t\in B_n,\psi(0)\in M\,and\,\psi|_{\partial B_n}=\psi_n|_{\partial B_n}.
\}
$$
Then for $j\geq2$, $c_j\geq c_*$, $K_{c_j}\setminus{W}\not=\emptyset$ and $c_j\rightarrow\infty$ as $j\rightarrow\infty$.
\end{theorem}

In order to apply Theorem \ref{Thm:xjl2}, we set $G = -id$, $J = I_{\lambda,\beta}$ and $P = P_\epsilon^+$. Then
$M = P_\epsilon^+
\cap P_\epsilon^-$, $\Sigma=\partial P_\epsilon^+
\cap \partial P_\epsilon^-$, and $W = P_\epsilon^+\cup P_\epsilon^-$.
 In this subsection, $f$ is
assumed to be odd, and, consequently, $ I_{\lambda,\beta}$ is even. Now, we show that $P_\epsilon^+$
is a G-admissible invariant set for the functional $ I_{\lambda,\beta}$ at any level $c$. Since $K_c$ is
compact, there exists a symmetric open neighborhood $N$ of $K_c \setminus{W}$ such that
$\gamma(\bar{N})<\infty$.

\begin{lemma}\label{Lem:deform1}
There exists $\epsilon_0 > 0$ such that for $0 <\epsilon<\epsilon'<\epsilon_0$, there exists a
continuous map $\sigma : [0,1]\times E \rightarrow E$ satisfying
\begin{itemize}
\item[\rm (1) ] $\sigma(0,u)=u$ for $u\in E$;
\item[\rm (2) ] $\sigma(t,u)=u$ for $t\in[0,1]$, $u\not\in I_{\lambda,\beta}^{-1}[c-\epsilon',c+\epsilon']$;
\item[\rm (3) ] $\sigma(t,-u)=-\sigma(t,u)$ for $(t,u)\in[0,1]\times E$;
\item[\rm (2) ]$\sigma(1,I_{\lambda,\beta}^{c+\epsilon}\setminus{(N\cup W)})\subset I_{\lambda,\beta}^{c-\epsilon}$;
\item[\rm (3) ] $\sigma(t,\overline{P_\epsilon^+})\subset\overline{P_\epsilon^+}$, $\sigma(t,\overline{P_\epsilon^-})
\subset\overline{P_\epsilon^-}$
for $t\in[0,1]$.
\end{itemize}
\end{lemma}
\Proof The proof is similar to that of Lemma \ref{Lem:XB}. Since $I_{\lambda,\beta}$ is even,
$B$ is odd
and thus $\sigma$ is odd in $u$.
\qed

\textbf{Proof of Theorem \ref{Thm:many} (Multiplicity)} We divide the proof into two steps.\\
\textbf{Step1.} Since $f$ is odd, it follows from Lemma \ref{Lem:deform1}
that $P_\epsilon^+$ is a G-admissible invariant set
for the functional $I_{\lambda,\beta}$ for $\lambda,\beta\in(0,1]$  at any level $c$.
We are now constructing $\psi_n$ satisfying the hypotheses of Theorem \ref{Thm:xjl2}.
 For any fixed $n\in\N$, we choose
$\{v_i\}_{i=1}^n\subset C_0^\infty(\R^3)\setminus\{0\}$ such that
$\text{supp}(v_i)\cap\text{supp}(v_j)=\emptyset$
for $i\not=j$. Define
$\psi_n\in C(B_n,E)$ as
$$
\psi_n(t)(\cdot)=R_n^2\Sigma_{i=1}^nt_iv_i(R_n\cdot),\quad t=(t_1,t_2,...,t_n)\in B_n.
$$
Observe that
$$
\rho_n=\min\{\|t_1v_1+t_2v_2+\cdot\cdot\cdot+t_nv_n\|_2:\,\,\sqrt{\Sigma_{i=1}^n t_i^2}= 1\}>0,
$$
then $\|u_t\|_2^2\geq\rho_n R_n$ for $u\in\varphi_n(\partial B_n)$ and
it follows from Lemma \ref{Lem:qfanshu} that
$\psi_i(\partial B_n)\cap P_\epsilon^+\cap P_\epsilon^-=\emptyset$. Similar to the proof of
Theorem \ref{Thm:existence} (existence part),
we also have
$$
\sup\limits_{u\in\psi_i(\partial B_n)}I_{\lambda,\beta}(u)<0<\inf\limits_{u\in\Sigma} I_{\lambda,\beta}(u).
$$
Clearly, $\psi_n(0)=0\in P_\epsilon^+\cap P_\epsilon^-$ and $\psi_n(-t)=-\psi_n(t)$ for $t\in B_n$.
For any fixed $\beta\in(0,1]$ and $j\in\{1,2,...,n\}$, we define
$$
c_{\lambda,\beta}^j=\inf\limits_{B\in\Gamma_j}\sup\limits_{u\in B\setminus{W}}I_{\lambda,\beta}(u),
$$
where $W:=P_\epsilon^+\cup P_\epsilon^-$ and $\Gamma_j$ was defined in Theorem \ref{Thm:xjl2}.
In view of the definition of $\Gamma_j$, it is easy to see that $c_{\lambda,\beta}^j$ is independent of $\epsilon$.
Based on Lemma \ref{Lem:xiaozhi} and Theorem \ref{Thm:xjl2}, for any fixed $\beta\in(0,1]$ and $j\geq2$,
\begin{equation}\label{eqn:fun9}
\frac{\epsilon^2}{4}\leq\inf\limits_{u\in\Sigma}I_{\lambda,\beta}(u):=c_*\leq c_{\lambda,\beta}^j\rightarrow\infty,\quad \text{as}\,\,j\rightarrow\infty
\end{equation}
and
there exists $\{u^j_{\lambda,\beta}\}\subset E\setminus{W}$ such that
$I_{\lambda,\beta}(u^j_{\lambda,\beta})=c_{\lambda,\beta}^j$ and $I'_{\lambda,\beta}(u^j_{\lambda,\beta})=0$.\\
\textbf{Step2.} Using  similar arguments to those in Theorem \ref{Thm:existence}, for any fixed
$j\geq 2$, $\{u^j_{\lambda,\beta}\}_{\lambda,\beta\in(0,1]}$ is bounded in $E$, that is to say,
there exists
$C>0$ independent of $\lambda,\beta$ such that $\|u^j_{\lambda,\beta}\|\leq C$. Without loss of generality,
we assume
$u^j_{\lambda,\beta}\rightharpoonup u^j_*$ weakly in $E$ as $\beta\rightarrow0^+$.
By Lemma \ref{Lem:xiaozhi} and Theorem \ref{Thm:xjl2} we have
$$
\frac{\epsilon^2}{4}\leq
\inf\limits_{u\in\Sigma} I_{\lambda,\beta}(u)\leq c_{\lambda,\beta}^j\leq c_{R_n}:=\sup\limits_{u\in\psi_n(B_n)}
I_{1,0}(u),
$$
where $c_{R_n}$ is independent of $\lambda,\beta$, and
$$
I_{1,0}(u):=\frac{1}{2}\int_{\R^3}(a|\nabla u|^2+V(x)u^2)+\frac{1}{2(1+\alpha)}\left(\int_{\R^3}|u|^2\right)^{1+\alpha}
-\int_{\R^3} F(u).
$$
Assume $c_{\lambda,\beta}^j\rightarrow c_*^j$
as $\lambda,\beta\rightarrow0^+$. Then we can prove
$u^j_{\lambda,\beta}\rightarrow u^j_*$ strongly in $E$ as $\lambda,\beta\rightarrow0^+$ and $u^j_*\in E\setminus{W}$
such that
$I'(u^j_*)=0$ and $I(u^j_*)=c_*^j$. We claim that $c_*^j\rightarrow\infty$ as $j\rightarrow\infty$. Indeed,
it follows from ($f_1$) and ($f_2$) that
\begin{equation}\label{eqn:modif-fanhan}
\aligned
I_{\lambda,\beta}(u)&\geq\frac{1}{2}\int_{\R^3}(a|\nabla u|^2+V(x)u^2)
-\int_{\R^3} F(u)-\frac{1}{r}\int_{\R^3}|u|^r\\
&\geq\frac{1}{2}\int_{\R^3}(a|\nabla u|^2+V(x)u^2)
-\int_{\R^3}(\frac{V_0}{4}u^2+\frac{C_{V_0}}{r}|u|^r)-\frac{1}{r}\int_{\R^3}|u|^r\\
&\geq\frac{1}{2}\int_{\R^3}(a|\nabla u|^2+W(x)u^2)
-\frac{C}{r}\int_{\R^3}|u|^r:=L(u),
\endaligned
\end{equation}
where $W(x):=V(x)-V_0/2$ and $C_{V_0},C>0$ are constants. We observe that
the boundedness of the Palais-Smale sequence is
not hard to verify for functional $L$ which satisfies
the corresponding Ambrosetti-Rabinowtiz condition. As a result,
with some suitable modification, the arguments of functional $I_{\lambda,\beta}$ are still
valid for $L$ without any perturbation. That is to say, the functional $L$ satisfies all conditions of Theorem \ref{Thm:xjl2}.
So we can define
\begin{equation}\label{eqn:modif-fanhan1}
d^j:=\inf\limits_{B\in\Gamma_j}\sup\limits_{u\in B\setminus{W}}L(u),
\end{equation}
where $W:=P_\epsilon^+\cup P_\epsilon^-$ and $\Gamma_j$ was defined in Theorem \ref{Thm:xjl2},
and $d^j$ is independent of $\epsilon$. And then, Theorem \ref{Thm:xjl2} gives $d^j\rightarrow+\infty$ as $j\rightarrow+\infty$.
Combining (\ref{eqn:modif-fanhan}) and (\ref{eqn:modif-fanhan1}),
it follows from
the definition of $c_{\lambda,\beta}^j>d^j$. Taking $\lambda,\beta\rightarrow0^+$,
we immediately get $c_*^j>b^j\rightarrow+\infty$ as $j\rightarrow+\infty$.
Therefore, problem (K) has infinitely many sign-changing solutions. The proof is complete.

\section{Energy doubling of sign-changing solutions}\setcounter{equation}{0}\label{sec5}
\renewcommand{\theequation}{5.\arabic{equation}}
In view of Theorem
\ref{Thm:existence}, we know that problem (K) has always a ground state sign-changing solution $w_b\in E$ for any $b>0$. We
prove now that $w_b$ is of an energy which is strictly large than $2c_b$ defined in (\ref{eqn:p-ground1}) as $b> 0$
is small.

\textbf{Proof of Theorem \ref{Thm:ground}}.
For any $b>0$, let $w_b\in E$ be a ground state sign-changing solution of
problem (K) with $I(w_b)=m_b$, where $m_b$ satisfies
\begin{equation}\label{eqn:mb}
m_b:=\inf\limits_{u\in\Theta}I(u),\quad\Theta:=\{u\in E\setminus\{0\},\,I'(u)=0,\,u^{\pm}\not\equiv0\}.
\end{equation}
In view of the proof of Theorem \ref{Thm:existence}, we can deduce from (\ref{eqn:fun0}) that
\begin{equation}\label{eqn:fun0--}
m_b\leq C_R:=\sup\limits_{u\in\psi_0(\triangle)}\bar{I}_{1}(u)<\infty,
\end{equation}
where $\psi_0(\triangle)$ was defined in the proof of Theorem \ref{Thm:existence} and
$C_R$ is independent of $b\in(0,1]$ and the functional $\bar{I}_{1}: E\rightarrow\R$ is defined as
$$
\bar{I}_{1}:=\frac{1}{2}\|u\|^2+\frac{1}{4}\left(\int_{\R^3}|\nabla u|^2\right)^2+\frac{1}{2(1+\alpha)}\|u\|_2^{2(1+\alpha)}-\int_{\R^3}F(u).
$$
We claim that, for any sequence $b_n\rightarrow0$, $\{w_{b_n}\}$ is a bounded sequence in $E$.
 We first have
\begin{equation}\label{eqn:fun1--}
m_{b_n}=\frac{a}{2}\int_{\R^3}|\nabla w_{b_n}|^2+\frac{1}{2}\int_{\R^3}V(x)w_{b_n}^2
+\frac{b_n}{4}\left(\int_{\R^3}|\nabla w_{b_n}|^2\right)^2-\int_{\R^3}F(w_{b_n})
\end{equation}
and
\begin{equation}\label{eqn:fun2--}
0=a\int_{\R^3}|\nabla w_{b_n}|^2+\int_{\R^3}V(x)w_{b_n}^2
+b_n\left(\int_{\R^3}|\nabla w_{b_n}|^2\right)^2-\int_{\R^3}f(w_{b_n})w_{b_n}.
\end{equation}
Moreover, from Lemma \ref{Lem:pohozaev}, the following identity holds
\begin{equation}\label{eqn:fun3--}
\aligned
\frac{a}{2}\int_{\R^3}|\nabla w_{b_n}|^2+&\frac{3}{2}\int_{\R^3}V(x)w_{b_n}^2+\frac{1}{2}\int_{\R^3}(\nabla
V(x),x)w_{b_n}^2\\+&\frac{b}{2}\left(\int_{\R^3}|\nabla w_{b_n}|^2\right)^2
-3\int_{\R^3}F(w_{b_n})=0.
\endaligned
\end{equation}
Multiplying (\ref{eqn:fun1--}), (\ref{eqn:fun2--}) and (\ref{eqn:fun3--})  by $4$, $-\frac{1}{\mu}$ and $-1$ respectively and adding them up,
we get
\begin{equation*}\label{eqn:fun4--}
\aligned
4m_{b_n}=&a\frac{3\mu-2}{2\mu}\int_{\R^3}|\nabla w_{b_n}|^2+\frac{\mu-2}{2\mu}\int_{\R^3}V(x)w_{b_n}^2
-\frac{1}{2}\int_{\R^3}(\nabla
V(x),x)w_{b_n}^2\\
&+\frac{\mu-2}{2\mu}b_n\left(\int_{\R^3}|\nabla w_{b_n}|^2\right)^2+\int_{\R^3}(\frac{1}{\mu}f(w_{b_n})w_{b_n}-F(w_{b_n})).
\endaligned
\end{equation*}
Combining with (\ref{eqn:fun0--}), we argue that $\{w_{b_n}\}$ is a bounded sequence in $E$.
Up to a subsequence, we assume that $w_{b_n}
\rightarrow w_0$ weakly in $E$. Note that $w_{b_n}$ is a ground state sign-changing solution
to problem (K) with $b=b_n$, then by the compactness of the embedding $E\hookrightarrow L^q(\R^3)$ $(2<q<6)$,
we deduce that $w_{b_n}
\rightarrow w_0$ strongly in $E$ and $w_0$ is a sign-changing solution
of (\ref{eqn:ellptic}).
Thus, by the energy doubling property of sign-changing solution of (\ref{eqn:ellptic}),
\begin{equation}\label{eqn:fun5--}
m_{b_n}=I(w_{b_n})=I_0(w_{0})+o(1)=m_0+o(1),
\end{equation}
where $m_0$ is given in (\ref{eqn:mb}) with $b=0$.

We now claim $m_0>2c_0$. Indeed, (\ref{eqn:p-ground-}) implies that $m_0\geq2c_0$.
Assume $m_0=2c_0$ and $w_0$ is a ground state sign-changing
solution of (\ref{eqn:ellptic}). Then we immediately get $w_0^{\pm}\in\mathcal{N}$ and
$I_0(w_0^{\pm})=c_0$. That is, $w_0^{\pm}$ is a minimizer of the functional $I_0$ which is restricted at
$\mathcal{N}$. It is known that $\mathcal{N}$ is a $C^1$ and natural constraint manifold. Hence,
by lagrange multiplier principle, we can easily get that
$w_0^{\pm}$ are critical points of the free functional $I_0$. This fact tells us that
$w_0^{+}$ and $w_0^{-}$ are two nonnegative solutions of (\ref{eqn:ellptic}).
Using the strong maximum principle, we can obtain two different positive solutions of (\ref{eqn:ellptic})
corresponding respectively to $w_0^{+}$ and $w_0^{-}$, contradicting to the uniqueness of positive solution.
Therefore, the claim is true. In view of (\ref{eqn:fun5--}), we have
\begin{equation}\label{eqn:fun5--}
m_{b_n}=I_0(w_{0})+o(1)=m_0+o(1)>2c_0,
\end{equation}
for large $n$.

We also claim that for each $b>0$,
there exists $u_b\in \mathcal{N}_b$ such that $I(u_b)=c_b$, and $u_b$ is a positive ground state solution of problem (K), where
$\mathcal{N}_b$ and $c_b$ have been given in (\ref{eqn:p-ground+}) and (\ref{eqn:p-ground1}), respectively.
In fact, similar to that of \cite{Li14,Liu15}, we construct a modified energy functional satisfying
the geometric conditions of monotonicity trick
developed by Struwe and Jeanjean \cite{Struwe85,Jeanjean99} to obtain a bounded
Palais-Smale sequence at a mountain-pass level. By the compactness of the embedding $E\hookrightarrow L^q(\R^3)$ $(2<q<6)$,
we can obtain a nontrivial critical point of the modified energy functional.
With the aid of the corresponding Pohozaev identity,
a convergence argument allows us to pass limit to the original problem (K) and
then to obtain a nontrivial solution in $E$. Without loss of generality,
we can assume that the nontrivial solution is
nonnegative, then the strong maximum principle implies that such a solution is positive.
We also obtain a positive ground state solution $u_b\in \mathcal{N}_b$ with $I(u_b)=c_b$,
when the above convergence argument is used to
the minimal sequence $\{u_n\}\subset\mathcal{N}_b$
satisfying $I(u_n)\rightarrow c_b$ as $n\rightarrow\infty$. The claim is proved.

Now taking $b_n\rightarrow0$, we
know that there exists $u_0\in E$ such that $u_{b_n}
\rightarrow u_0$ strongly in $E$ and $u_0$ is a
positive solution of problem (\ref{eqn:ellptic}).
Hence, by the uniqueness of positive solution of (\ref{eqn:ellptic}), we have
\begin{equation}\label{eqn:fun6--}
c_{b_n}=I(u_{b_n})=I_0(u_0)+o(1)=c_0+o(1),
\end{equation}
where $c_0$ is given in (\ref{eqn:p-ground}). Combining (\ref{eqn:fun5--}) with (\ref{eqn:fun6--}), we have
$m_{b_n}>2c_{b_n}$ for large $n\in\N$. Hence there exists $b^*>0$ such that $m_{b}>2c_{b}$ for any $b\in(0,b^*)$.
By the definition of $c_{b}$, $m_{b}$ is strictly two times larger than that of the ground state energy.
The proof is complete.
\qed

\end{document}